%% file: piercing.tex
\documentclass[11pt]{amsart}
\usepackage{amsmath,amssymb, amsthm, mathabx,amscd, amsfonts, amssymb, hyperref, mathrsfs, textcomp, epsfig, graphicx, IEEEtrantools, tikz, verbatim, xypic, subcaption, a4wide, arydshln}
\usepackage[mathscr]{euscript}

\usetikzlibrary{shapes,arrows,cd}
\usetikzlibrary{matrix,decorations.pathreplacing,angles,quotes}
\usetikzlibrary{decorations.markings}

 \tikzset{->-/.style={decoration={
  markings,
  mark=at position .5 with {\arrow{>}}},postaction={decorate}}}

\newtheorem{thm}{Theorem}
\newtheorem{prop}{Proposition}[section]
\newtheorem{lemma}{Lemma}[section]

\newtheorem*{quest*}{Main Question}

\newcommand{\spin}{\mathfrak{s}}
\theoremstyle{definition}
\newtheorem{defn}{Definition}[section]

\newtheorem{cond}{Condition}

\theoremstyle{remark}
\newtheorem{remark}{Remark}[section]
\newtheorem{example}{Example}[section]

     \RequirePackage{rotating}                   
    \def\HSt{%
       \setbox0=\hbox{$\widehat{\mathit{HS}}$}
       \setbox1=\hbox{$\mathit{HS}$}
       \dimen0=1.1\ht0
       \advance\dimen0 by 1.17\ht1
       \smash{\mskip2mu\raise\dimen0\rlap{%
          \begin{turn}{180}
              {$\widehat{\phantom{\mathit{HS}}}$}
           \end{turn}} \mskip-2mu    
                \mathit{HS}
    }{\vphantom{\widehat{\mathit{HS}}}}{}}

     \RequirePackage{rotating}                   
    \def\HMt{%
       \setbox0=\hbox{$\widehat{\mathit{HM}}$}
       \setbox1=\hbox{$\mathit{HM}$}
       \dimen0=1.1\ht0
       \advance\dimen0 by 1.17\ht1
       \smash{\mskip2mu\raise\dimen0\rlap{%
          \begin{turn}{180}
              {$\widehat{\phantom{\mathit{HM}}}$}
           \end{turn}} \mskip-2mu    
                \mathit{HM}
    }{\vphantom{\widehat{\mathit{HM}}}}{}}

\newcommand{\HMred}{{\mathit{HM}}}

\newcommand{\sinc}{\mathrm{sinc}}
\newcommand{\vol}{\mathrm{vol}}

\newcommand{\T}{\mathcal{T}_+}

\newcommand{\hol}{\mathrm{hol}}

\begin{document}

\title[Dirac spectra and Floer homology]{Dirac spectral flow and Floer theory of hyperbolic three-manifolds}

\author{Francesco Lin}
\address{Department of Mathematics, Columbia University} 
\email{flin@math.columbia.edu}

\author{Michael Lipnowski}
\address{Department of Mathematics, Ohio State University} 
\email{michaellipnowski@gmail.com}

\begin{abstract}
We study the interplay between hyperbolic geometry and monopole Floer homology for a closed oriented three-manifold $Y$ with $b_1=1$ equipped with a torsion spin$^c$ structure $\spin$. We show that, under favorable circumstances, one can completely describe the Floer theory of $(Y,\spin)$ purely in terms of geometric data such as the lengths and holonomies of closed geodesics. In particular, we perform the first computations of monopole Floer chain complexes with non-trivial homology for hyperbolic three-manifolds.
\par
The examples we consider admit no irreducible solutions to the Seiberg-Witten equations, and the non-triviality of the Floer homology groups is a consequence of the geometry of the $1$-parameter family of Dirac operators associated to flat spin$^c$ connections. The main technical challenge is to understand explicitly how the Dirac eigenvalues with small absolute value cross the value zero in this family; we tackle this using Fourier analytic tools via the corresponding $1$-parameter family of odd Selberg trace formulas and its derivative.
\end{abstract}

\maketitle

\section*{Introduction}
A fundamental outstanding problem in low-dimensional topology is to understand the interactions between hyperbolic geometry and Floer theory. A concrete (open ended) instance of this is the following.
\begin{quest*}
For a given closed oriented hyperbolic three-manifold $Y$, can one determine the Floer homology groups of $Y$ in terms of  the numerical invariants arising from hyperbolic geometry (such as volume, lengths of closed geodesic, etc.)?
\end{quest*}
Throughout the paper we will focus on the monopole Floer homology package \cite{KM}, to which we refer simply as `Floer homology'. In \cite{LL1}, the authors described a method that can be applied in favorable situations to show that a given hyperbolic rational homology sphere $Y$ is an $L$-space, i.e. its Floer homology is trivial. This is achieved by using the spectral theory of $Y$, and more specifically the spectral gap $\lambda_1^*$ of the Hodge Laplacian acting on coexact $1$-forms, as a stepping stone between Floer theory and hyperbolic geometry. On one hand, it is shown that if $\lambda_1^*>2$ (i.e. if $Y$ is \textit{spectrally large}), then the Seiberg-Witten equations on $Y$ do not admit irreducible solutions. On the other hand, the following version of the Selberg trace formula
\begin{align}\label{tracecoexact}
\begin{split}
\overbrace{\frac{1}{2}(b_1(Y)-1) \cdot \widehat{H}(0)+\frac{1}{2}\sum_{j\geq 1} \widehat{H}(\sqrt{\lambda_j^*})}^{\text{spectral side}}=&\frac{\mathrm{vol(Y)}}{2\pi} \cdot (H(0)-H''(0))\\
&\underbrace{+\sum_{\gamma} \ell(\gamma_0)\frac{\cos(\mathrm{hol}(\gamma))}{|1-e^{\mathbb{C}\ell(\gamma)}||1-e^{-\mathbb{C}\ell(\gamma)}|}H(\ell(\gamma))}_{\text{geometric side}}
\end{split}
\end{align}
holds for any even, compactly supported, sufficiently smooth test function $H$, where:
\begin{itemize}
\item $\widehat{H}(t)$ denotes the Fourier transform of $H(x)$;
\item the sum on the spectral side is taken over the spectrum
\begin{equation*}
0<\lambda_1^*\leq \lambda_2^*\leq\dots
\end{equation*}
of the Hodge Laplacian $\Delta=(d+d^*)^2$ acting on coexact $1$-forms;
\item  the sum on the geometric side is taken over closed geodesics $\gamma$ and \begin{equation}\label{complexlength}
\mathbb{C}\ell(\gamma)=\ell(\gamma)+i\mathrm{hol(\gamma)\in \mathbb{R}+i(\mathbb{R}/2\pi\mathbb{Z})}
\end{equation}
denotes the complex length, while $\gamma_0$ a prime geodesic $\gamma$ is multiple of.
\end{itemize}
Adapting to hyperbolic three-manifolds the Fourier optimization techniques from \cite{BS}, used there in the context of the trace formula on hyperbolic surfaces (with a view towards the celebrated Selberg-$1/4$ conjecture), this formula can be used to provide explicit lower bounds on $\lambda_1^*$ in terms of the volume of $Y$ and the complex length spectrum computed (using SnapPy \cite{SnapPy}, for example) up to a given cutoff $R>0$. The method is successful in proving that many small volume manifolds in the Hodgson-Weeks census \cite{Census} are $L$-spaces; see also \cite{LL2} for a larger volume example.
\\
\par
A significant limitation of this approach is that it can only be used to show that a given manifold has trivial Floer homology. In the present work, we will instead provide the first computations of Floer chain complexes of spin$^c$ hyperbolic three-manifolds $(Y,\spin)$ whose homology is \textit{not} trivial. While the Floer homology of our examples can be also computed using more standard topological tools, the key point is that our approach takes as input \textit{solely} numerical quantities arising from hyperbolic geometry, in the spirit of the Main Question.
\par
Our examples consist of hyperbolic three-manifolds $Y$ with $b_1=1$ which are spectrally large in the same sense as above, i.e. $\lambda_1^*>2.$  While \cite{LL1} focuses on the case $b_1=0$, the same approach readily applies to certify that $\lambda_1^\ast > 2$ on manifolds for which $b_1>0$ too.  We equip $Y$ with a spin$^c$ structure $\spin$ which is torsion, i.e. $c_1(\spin)$ is torsion. Recall that such structures are (not canonically) in bijection with the torsion subgroup of $$H^2(Y;\mathbb{Z})\stackrel{PD}{\equiv}H_1(Y;\mathbb{Z}).$$
For simplicity, in this introduction we focus on the Floer homology groups $\HMred_*(Y,\spin;\Gamma_\eta)$ with coefficients in the $\mathbb{R}$-local system $\Gamma_\eta$ on the moduli space of configurations corresponding to a real $1$-cycle $\eta\subset Y$ with $$[\eta]\neq 0\in H_1(Y;\mathbb{R}),$$ see \cite[Section 3.7]{KM}. This is an absolutely $\mathbb{Q}$-graded module over $\mathbb{R}[U]$, with $U$ acting with degree $-2$, finite dimensional as an $\mathbb{R}$-vector space. Our techniques in fact determine the `usual' Floer homology $\HMt_*(Y,\spin)$ as well, cf. Section \ref{background} for the corresponding expressions; in fact, we will compute explicitly the Floer chain complexes for these manifolds, both with twisted and untwisted coefficients. Focusing again on manifolds in the Hodgson-Weeks census, we obtain the following examples in which $\HMred_*(Y,\spin;\Gamma_\eta)$ is either $0$ or $\mathbb{R}$.

\begin{thm}\label{thm1}We can determine explicitly the monopole Floer chain complexes for the torsion spin$^c$ structures on the census manifolds satisfying $b_1(Y)=1$ which are listed in Table 1 (all of which are spectrally large). The corresponding Floer homology groups with local coefficients are either$$\HMred_*(Y,\spin;\Gamma_\eta)=0\text{ or }\mathbb{R},$$ the number of occurrences of the latter case (among torsion spin$^c$ structures) being indicated in the last column.
\end{thm}

\begin{table}[ht]
\centering
\begin{tabular}{|c| c |c |c|c|}
\hline
Census label & volume & systole &  $H_1(Y;\mathbb{Z})_{\mathrm{tors}}$&  $\#\mathbb{R}$ \\ [0.5ex] 
\hline
356&3.1663&0.3887&$\mathbb{Z}/3$& 1\\
357&3.1663&0.3887&$\mathbb{Z}/7$& 3\\
381&3.1772&0.3046&$\mathbb{Z}/7$&4\\
734 & 3.6638&1.0612&$\mathbb{Z}/2$&0\\
735&3.6638&0.5306&$\mathbb{Z}/10$& 4\\
790&3.7028&0.5585& \{0\} &0\\
882&3.7708&0.3648&$\mathbb{Z}/11$& 5\\
1155 & 3.9702 & 0.3195 &$\mathbb{Z}/6$ &4\\
1280&4.0597&0.5435&$\mathbb{Z}/13$&5\\
1284& 4.0597& 0.4313&$\mathbb{Z}/4$&1\\
3250&4.7494&0.5577&$\mathbb{Z}/2$&1\\
3673& 4.8511& 0.3555 &$\mathbb{Z}/17$&8 \\
\hline
\end{tabular}
\label{examplesthm}
\caption{The spectrally large manifolds in Theorem \ref{thm1}, with approximate volumes and systoles (i.e. length of shortest closed geodesic).}
\end{table}

Among the $\approx11k$ manifolds in the Hodgson-Weeks census, only $127$ satisfy $b_1=1$, and we certified that at least $\approx 30$ of them are spectrally large. While the strategy we will employ in the proof of Theorem \ref{thm1} can be applied to determine the Floer homology of many other spectrally large examples (provided enough length spectrum is computed), for concreteness we focused our attention on the ten examples with smallest volume (having census labels between 356 and 1284), all of which are fibered with Thurston norm $2$, and the two smallest non-fibered examples (having census labels 3250 and 3673), which also have Thurston norm $2$, cf. Table $7.7$ in \cite{Chen}.
\par
In fact, our techniques can also be applied to compute Floer chain complexes having higher rank homology. As a concrete example of this we will prove the following.

\begin{thm}\label{thm2}
The census manifold \#10867, which is spectrally large and has $H_{1}(Y;\mathbb{Z})=\mathbb{Z}\oplus\mathbb{Z}/4$, admits a self-conjugate spin$^c$ structure $\spin$ for which we can determine explicitly the monopole Floer chain complex and for which $$\HMred_*(Y,\spin;\Gamma_\eta)=\mathbb{R}^{\oplus 2},$$ with both summands lying in the same absolute grading. For the other self-conjugate spin$^c$ structure, the Floer homology is $\mathbb{R}$, while for the non-self-conjugate ones it is trivial.
\end{thm}
Here, self-conjugate means that $\bar{\spin}=\spin$, or equivalently that $\spin$ is induced by a spin structure; in fact, $\spin$ is induced by exactly two spin structures because $b_1=1$. The manifold \#10867, which is fibered with Thurston norm $4$ (cf. Table $7.7$ in \cite{Chen}), has volume $\approx 6.2391$ and systole $\approx 0.3457$.

\begin{remark}
Using the techniques developed in \cite{LL3}, one can also determine the absolute $\mathbb{Q}$-grading \cite[Ch. 28]{KM} in terms of numerical quantities arising in hyperbolic geometry; we will not pursue such computations in this paper.
\par
The Floer homology groups also admit a canonical $\mathbb{Z}/2\mathbb{Z}$-grading, and all the groups in Theorems \ref{thm1} and \ref{thm2} lie in even grading. Recall also that the Euler characteristic of these groups with respect to this grading is the same as the Turaev torsion of $(Y,\spin)$ \cite{Tur}.
\end{remark}

\begin{remark}
A natural source of manifolds with $b_1=1$ are longitudinal surgeries on knots in $S^3$.  Unfortunately, we did not find spectrally large examples among the hyperbolic ones obtained by surgery on small crossing knots.
\end{remark}

The proof of Theorems \ref{thm1} and \ref{thm2} again involves spectral theory as a link between Floer theory and hyperbolic geometry. Unlike the case of $b_1=0$, the outcome strongly depends on the choice of spin$^c$ structure, and an essential role is played by the spectra of Dirac operators and how they vary in a $1$-parameter family. We will now give an overview of how the interactions play out when studying the problem under consideration.
\\
\par
\begin{figure}
\centering
  \includegraphics[width=.8\linewidth]{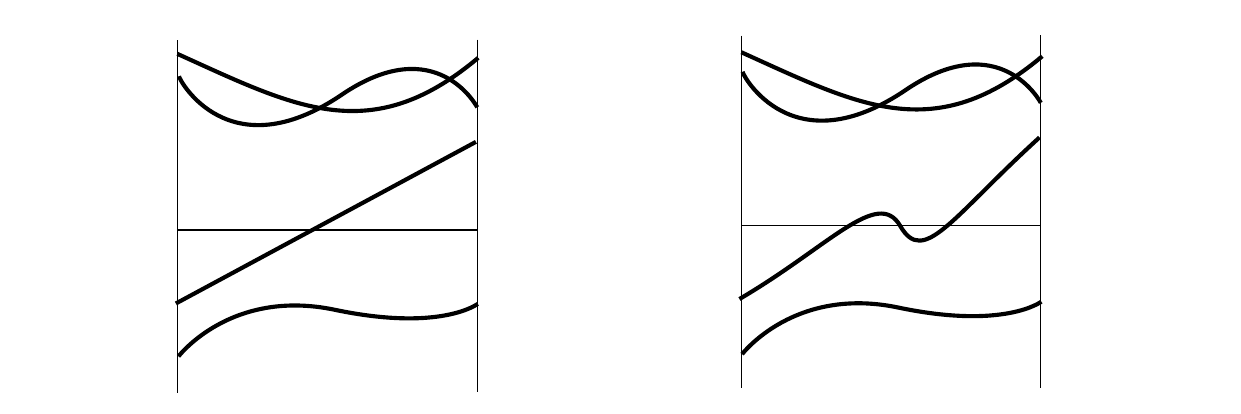}
  \caption{On the left, the small eigenvalue (i.e. the eigenvalue with smallest absolute value) crosses zero with sign $+$, while on the right it crosses it three times with signs $+,-,+$.}
  \label{spectralflowpic}
\end{figure}

\textbf{From spectral theory to Floer homology.} Under the hypothesis that $Y$ is spectrally large, the same argument as in \cite{LL1} shows that the unperturbed Seiberg-Witten equations do not admit irreducible solutions. On the other hand, the reducible solutions - consisting of flat spin$^c$ connections $B$ up to gauge - are parametrized by 
\begin{equation}\label{torus}
\mathbb{T}={H^1(Y;i\mathbb{R})}/{H^1(Y;2\pi i\mathbb{Z})}, 
\end{equation}
which is a circle because $b_1 = 1$, and perturbing the equations to achieve transversality is thus significantly subtler. In this sense, a key role is played by the family of Dirac operators $$\{D_B\} \text{ with }[B]\in\mathbb{T},$$ each of which is a self-adjoint operator diagonalizable in $L^2$ with discrete real spectrum unbounded in both directions. More specifically, we will be concerned with the geometry of the locus
\begin{equation}\label{K}
\mathsf{K}=\{[B]\text{ such that }\ker D_B\neq\{0\}\}\subset\mathbb{T}
\end{equation}
consisting of operators with non-trivial kernel. Generically, one expects this to be a collection of points with attached signs $\pm$ that record whether the eigenvalue of smallest absolute value goes from negative to positive or vice versa; the latter is usually referred to as the spectral flow at the crossing \cite{APS3}, see Figure \ref{spectralflowpic}. The significance of $\mathsf{K}$ is that it is exactly the locus in $\mathbb{T}$ at which the Chern-Simons-Dirac functional fails to be Morse-Bott.

\begin{remark}
Throughout the text, we will refer to an eigenvalue with small absolute value as a \textit{small} eigenvalue.
\end{remark}

We will refer to the sequence of signs that one encounters going around $\mathbb{T}$ as the \textit{piercing sequence} of $(Y,\spin)$. It turns out that in the spectrally large case, the piercing sequence of $(Y,\spin)$ completely determines $\HMred_*(Y,\spin;\Gamma_\eta)$: in Theorems \ref{thm1} and \ref{thm2}, the Floer homology is $0$, $\mathbb{R}$ or $\mathbb{R}^{\oplus2}$ corresponding to the piercing sequence being empty, $(+,-)$, or $(+,-,+,-)$ respectively. Indeed, one can describe explicitly the corresponding Floer chain complexes, see Section \ref{background}.

\begin{remark}
In the situation we have just described, the non-triviality of the Floer homology of $(Y,\spin)$ is a consequence of the fact that the spectrum of the Dirac operator $D_B$ depends on the specific connection $[B]\in\mathbb{T}$, rather than the curvature $F_{B^t}$, which is identically zero. This is exactly the mechanism behind the Aharonov-Bohm effect in physics \cite{SN}.
\end{remark}

\vspace{0.3cm}
\textbf{From hyperbolic geometry to spectral theory.} Most of the work in the
 present paper is devoted to developing techniques to determine the geometry of the locus $\mathsf{K}$ and in particular the piercing sequence of $(Y,\spin)$. The key tool we will use is a version of the the Selberg trace formula for the Dirac operator $D_B$ associated with a flat spin$^c$ connection $B$ proved in \cite{LL3}. The formula, which has the same overall structure of (\ref{tracecoexact}), provides for each test function (either even or odd) an identity relating the spectrum $\{s_i\}_{i\in\mathbb{Z}}$ of $D_B$ to the geometry of $(Y,\spin)$ in the form of the \textit{spin$^c$ length spectrum} of $Y$ corresponding to the connection $B$. The latter is a refinement of the complex length spectrum (\ref{complexlength}) which takes into account the specific flat spin$^c$ connection $B$, and can be computed in concrete examples starting from a Dirichlet domain of $Y$ (as provided from SnapPy) using the techniques developed in \cite{LL3}.
\par
In fact, the techniques of our previous papers can be adapted to determine a collection of (small) disjoint, closed intervals $\{I_k\}$ in $\mathbb{T}$ (parametrized by $\tau$) for which:
\begin{itemize}
\item[(a)]
if $B_\tau \in \mathsf{K},$ then necessarily $\tau \in \bigsqcup I_k,$ and
\item[(b)]
the Dirac operators corresponding to $\tau \in \bigsqcup I_k$ admit a unique small eigenvalue, which we denote by $s_0(\tau).$
\end{itemize}
Furthermore, we find slightly larger closed disjoint intervals $\tilde{I}_k \supset I_k$ such that $s_0(\tau)$ becomes large and of opposite signs at the endpoints of $\tilde{I}_k.$  By continuity, $s_0(\tau_0)=0$ for some $\tau_0\in I_k.$
\par
The new challenge is to prove that, in each interval $I_k$, the smallest eigenvalue crosses zero \textit{at exactly one point, in a transverse fashion}, see Figure \ref{spectralflowpic}. This is important for our purposes because the Floer homology depends on the crossings via the piercing sequence; it is challenging because the techniques of our previous papers only provide $C^0$-information about $s_0(\tau)$, while we need to obtain information about the derivative $s_0'(\tau)$ in the present context. In particular, it would suffice to show that
\begin{equation}\label{nonzeroder}
s_0'(\tau)\neq 0 \text{ for } \tau\in I_k.
\end{equation}
Notice that if this holds, the sign attached to $D_{B_{\tau_0}} \in \mathsf{K}$ equals the constant sign of $s_0'(\tau)$ for all $\tau \in I_k$.
\par
A key observation is that the $1$-parameter family of trace formulas for $\{D_B\}_{[B]\in\mathbb{T}}$ can also be used to obtain information about the derivatives of the eigenvalues in the family. Roughly speaking, our strategy to prove (\ref{nonzeroder}) involves, for a fixed odd test function $K$, considering the \textit{derivative in $\tau$} of the $1$-parameter family of trace formulas associated with $\{D_B\}_{[B]\in\mathbb{T}}$. From this, we obtain a relation between all eigenvalues $\{s_i(\tau)\}$ of $D_B$, their derivatives $\{s_i'(\tau)\}$ and the hyperbolic geometry of $Y$.
\par
When proving (\ref{nonzeroder}), most of the technical work will be devoted to bounding from above the contribution coming from the eigenvalues $s_i$ with $i\neq 0$. A fundamental problem in the latter process is to obtain an a priori \textit{explicit} upper bound on the derivatives of the eigenvalues $s_i'$. As we will explain, such a bound can be obtained in terms of the $C^0$-norm of any closed $1$-form representing the generator of $H^1(Y;\mathbb{Z})=\mathbb{Z}$; this is the last numerical input from hyperbolic geometry required in our proofs, and we discuss in Appendix \ref{C0bounds} how to determine such bounds explicitly in our examples of interest. 
\\
\par
\textbf{Computational aspects.} In principle, our approach works for any spectrally large manifold with $b_1=1$ (given of course the statement to be proved is true) provided the computation of the length spectrum up to a given cutoff $R$ large enough. Of course, this is infeasible for large $R$ because of the prime geodesic theorem
\begin{equation*}
\#\{\text{prime geodesics $\gamma$ with $\ell(\gamma)\leq R$} \}\approx \frac{e^{2R}}{2R}\text{ for }R\rightarrow\infty,
\end{equation*}
see for example \cite{Cha}. We will focus on examples for which the computation of the length spectrum up to $R\in [7,8.5]$ is enough to implement our strategy; depending on the specific situation, the computation of the relevant spectrum took somewhere between a few hours and ten days of CPU time. The key feature that makes our approach work is that in the intervals with small eigenvalues, there are no other eigenvalues with small absolute value, or, if there are, their number and position can be understood somewhat explicitly.
\\
\par
\textit{Organization of the paper.} In Section \ref{background} we discuss the relevant background in monopole Floer homology, while in Section \ref{traceformulas} we recall the trace formulas for Dirac operators needed for our purposes, and derive the formula for the derivatives of the eigenvalues we will use.
\par
In the following two sections, we discuss in detail the techniques involved in the proof of our main theorems, working them out in detail for the unique self-conjugate spin$^c$ structure on the Census manifold \#357. In particular, in Section \ref{smalleigint} we build on our previous work to find the intervals containing small eigenvalues, while in Section \ref{singlecrossing} we prove that in each of these intervals there is exactly a single transverse crossing. In general, the strategy needs to be tweaked to be applied to different examples, and in Section \ref{moreexamples} we discuss in detail how to overcome some difficulties in the examples \#3250 and \#10867.
\par
Finally, in Appendix \ref{C0bounds} we discuss how to construct explicitly closed $1$-forms corresponding to generators of $H^1(Y;\mathbb{Z})$ for our manifolds of interest, and bound their $C^0$ norm. 
\\
\par
\textit{Acknowledgements.} The first author was partially supported by NSF grant DMS-2203498.  The second author was partially supported by NSF CAREER grant 2338933.
\vspace{0.5cm}

\section{The Floer theoretic setup}\label{background}
In this section we recollect the relevant background regarding monopole Floer homology; the fundamental reference for the topic is \cite{KM}, see also \cite{LinLec} for a friendly introduction. We will consider a Riemannian three-manifold $(Y,g)$ equipped with a torsion spin$^c$ structure $\spin$. We will consider the Floer homology group $\HMt_*(Y,\spin)$ and its version $\HMred(Y,\spin;\Gamma_{\eta})$ corresponding to the $\mathbb{R}$-local coefficient system $\Gamma_\eta$ associated with a $1$-cycle $\eta$ with $[\eta]\neq0$, see \cite[Section $3.7$]{KM}.
\\
\par
The unperturbed Seiberg-Witten equations for a pair $(B,\Psi)$ consisting of a spin$^c$ connection $B$ and a spinor $\Psi$ are
\begin{align}\label{seibergwitten}
\begin{split}
D_B\Psi&=0\\
\frac{1}{2}\rho(\ast F_{B^t})+(\Psi\Psi^*)_0&=0.
\end{split}
\end{align}
Here $D_B$ is the Dirac operator corresponding to $B$, $\rho$ is the Clifford multiplication and $F_{B^t}$ is the curvature of the connection $B^t$ induced on the determinant line bundle of the spinor bundle $S$. A solution $(B,\Psi)$ is called reducible if $\Psi\equiv 0$, and irreducible otherwise. By Hodge theory, reducible solutions up to gauge can be identified with the $b_1(Y)$-dimensional torus $\mathbb{T}$ of flat spin$^c$ connections (\ref{torus}).
\par
The monopole Floer homology groups of $(Y,\spin)$ are obtained by suitably counting solutions to (\ref{seibergwitten}) up to gauge, in the sense of the $S^1$-equivariant homology of the Chern-Simons-Dirac functional $\mathcal{L}$ (where the $S^1$, thought of as the group of constant gauge transformations acting on the spinor by multiplication). For example, the version $\HMt_*(Y,\spin)$ corresponds to Borel homology.
\\
\par
A key observation for our purposes is that while the locus $\mathbb{T}$ of reducible solutions is a smooth manifold, it fails in general to be non-degenerate in a Morse-Bott sense; indeed, the Hessian of $\mathcal{L}$ along $\mathbb{T}$ is singular in the normal directions exactly along the locus (\ref{K}) at which the corresponding Dirac operator has kernel. Therefore, in order to achieve transversality in the sense of \cite{KM}, one needs to appropriately perturb the equations so that there are only finitely many reducible solutions $[B,0]$, with the further property that $D_B$ has simple spectrum and no kernel. For the definition of the version $\HMt_*(Y,\spin)$, each such (perturbed) reducible solution contributes to the Floer chain complex a so-called tower
\begin{equation*}
\T=\mathbb{Z}[U,U^{-1}]/\mathbb{Z}[U],
\end{equation*}
where the copy of $\mathbb{Z}$ labeled $U^{-i}$ is generated by the stable critical point corresponding to the $i$th positive eigenvalue of $D_B$.
\par
We will be interested in the case in which the Seiberg-Witten equations on $(Y,\spin)$ have no irreducible solutions. When $b_1=0$ this implies (assuming that the Dirac operator at the unique reducible solution has no kernel) that we simply have $$\HMt_*(Y,\spin)=\T$$ as $\mathbb{Z}[U]$-modules. On the other hand, it turns out that when $b_1>0$, the Floer homology of $(Y,\spin)$ might be non-trivial even in the absence of irreducible solutions; indeed, it will crucially depend on the geometry of the locus $\mathsf{K}\subset\mathbb{T}$. In what follows, let us focus on the case of a hyperbolic manifold $Y$ with $b_1(Y)=1$, which is the one relevant for our purposes. We will consider two conditions.
\begin{cond}The first eigenvalue $\lambda_1^*$ of the Hodge Laplacian $\Delta=(d+d^*)^2$ acting on coexact $1$-forms on $Y$ satisfies $\lambda_1^*>2$.
\end{cond}
We will refer to this condition as being \textit{spectrally large}, and it can be verified to hold for a given hyperbolic three-manifold $Y$ via the work \cite{LL1} involving the Selberg trace formula for coexact $1$-forms (\ref{tracecoexact}). Here the threshold arises from the fact that hyperbolic three-manifolds have constant Ricci curvature $-2$.
\par
While the condition of being spectrally large is independent of the spin$^c$ structure, the following one is not.
\begin{cond}
The locus $\mathsf{K}\subset \mathbb{T}$ is transversely cut out in the sense that it consists of a finite number of points corresponding to operators $D_{B_{\tau_0}}$ for $\tau_0\in \mathbb{T}$ with simple kernel such that, denoting by $s:I\rightarrow \mathbb{R}$ a parametrization of the eigenvalue with smallest absolute value of $D_{B_\tau}$ in a small interval $I$ containing $\tau_0$, $s$ is smooth and transverse to $0$.
\end{cond}
In particular, one can assign to each point in $\mathsf{K}$ a $\pm$ sign depending of whether the eigenvalue of smallest absolute value goes from negative to positive or viceversa.
\begin{defn}
If Condition $2$ holds, we define the \textit{piercing sequence} of $(Y,\spin)$ to be the sequence of signs one encounters going around $\mathbb{T}$, up to cyclic reordering.
\end{defn}
\begin{remark}
Given that $\spin$ is torsion, the sum of plus and minuses along the circle is zero by the Atiyah-Singer index theorem for families \cite{APS3}.
\end{remark}
In \cite{LinDir} it is shown how to suitably perturb the Seiberg-Witten equations to achieve transversality while not introducing irreducible solutions when Conditions $1$ and $2$ hold. Roughly speaking, this involves a Morse function $f:\mathbb{T}\rightarrow \mathbb{R}$ achieving the same value at all maxima and minima, and such that at the points of $\mathsf{K}$ the vector field $\mathrm{grad} f$ is the direction in which the eigenvalue of smallest eigenvalue increases. The reducible critical points of the perturbed equations correspond to the critical points of $f$ on $\mathbb{T}$.
\par
A consequence of this is that one can compute the Floer chain complex in a completely explicit fashion; the key phenomenon that gives rise to interesting Floer homology is that the presence of points of $\mathsf{K}$ causes the various towers associated to reducible critical points to shift in grading because of spectral flow.
\begin{remark}
Even in situations where the chain complex $\check{C}_*$ only involves (stable) reducible critical points, in general the Floer differential
\begin{equation*}
\check{\partial}=\bar{\partial}^s_s-\partial^u_s\bar{\partial}^s_u
\end{equation*}
also involves irreducible trajectories via the term $\partial^u_s$. On the other hand, when $b_1=1$ the second term vanishes for grading reasons, so the $\check{\partial}$ only involves reducible trajectories.
\end{remark}
\subsection{Concrete examples.}\label{exampleschain} Let us describe in detail the output in the situations relevant for our purposes, namely those for which the piercing sequences are empty, $(+,-)$ and $(+,-,+,-)$. A more general discussion for self-conjugate spin$^c$ structures, relevant in the case of three-manifolds with $b_1=1$ obtained by mapping tori of finite order mapping classes, can be found in \cite{LinTheta}.
\par
In what follows, we will denote by $M\langle d\rangle$ the module obtained from $M$ by shifting degree upwards by $d$, i.e. $M\langle d\rangle_n=M_{n-d}$; we will also describe the action of the generator $\gamma\in H_1(Y;\mathbb{Z})/\mathrm{tors}$ (which has degree $-1$).

\begin{example}If the piercing sequence is empty, we have
\begin{equation*}
\check{C}_*(Y,\spin)=\T\oplus \T\langle1\rangle
\end{equation*}
with trivial differential, so that
\begin{equation*}
\HMt_*(Y,\spin)=\T\oplus \T\langle1\rangle
\end{equation*}
and the action of $\gamma$ is an isomorphism from the right tower onto the left tower. Correspondingly, we have for a local coefficient system $\Gamma_\eta$ with $[\eta]$ non-zero
\begin{equation*}
\HMred_*(Y,\spin;\Gamma_\eta)=0.
\end{equation*}
Note that this is the same as the case of $(S^1\times S^2,\spin_0)$ discussed in \cite[Ch. 36]{KM}.
\end{example}

\begin{example}If the piercing sequence is $(+,-)$, we have
\begin{equation*}
\check{C}_*(Y,\spin)=\T\oplus \T\langle-1\rangle
\end{equation*}
with trivial differential, so that
\begin{equation*}
\HMt_*(Y,\spin)=\T\oplus \T\langle-1\rangle
\end{equation*}
and the action of $\gamma$ is surjective from the right tower onto the left tower, and zero on the bottom of the right tower. The downward shift of the right tower comes from spectral flow. Correspondingly, we have for a local coefficient system $\Gamma_\eta$ with $[\eta]$ non-zero
\begin{equation}\label{nontrivial}
\HMred_*(Y,\spin;\Gamma_\eta)=\mathbb{R}_{-1}.
\end{equation}
Note that this is the same as the case of the `interesting' spin$^c$ structure on flat manifolds with $b_1=1$ discussed in \cite[Ch. 37]{KM}.
\end{example}

\begin{remark}\label{chainpicture}More pictorially, the previous two examples (using $\mathbb{R}$ coefficients for the `usual' untwisted Floer homology as well) are represented by the following diagrams respectively

\begin{center}
\begin{tikzcd}[row sep=small]
	& \vdots &&&&&& \vdots \\
	\vdots & {\mathbb{R}} &&&&& \vdots & {\mathbb{R}} \\
	{\mathbb{R}} &&&&&& {\mathbb{R}} \\
	& {\mathbb{R}} &&&&&& {\mathbb{R}} \\
	{\mathbb{R}} &&&&&& {\mathbb{R}} \\
	& {\mathbb{R}} &&&&&& {\mathbb{R}} \\
	{\mathbb{R}} &&&&&& {\mathbb{R}} \\
	&&&&&&& {\underline{\mathbb{R}}}
	\arrow["{1-e^{[\eta]}}", from=2-2, to=3-1]
	\arrow["{1-e^{[\eta]}}", from=2-8, to=3-7]
	\arrow["{1-e^{[\eta]}}", from=4-2, to=5-1]
	\arrow["{1-e^{[\eta]}}", from=4-8, to=5-7]
	\arrow["{1-e^{[\eta]}}", from=6-2, to=7-1]
	\arrow["{1-e^{[\eta]}}", from=6-8, to=7-7]
\end{tikzcd}
\end{center}
where we consider $[\eta]\in H_1(Y;\mathbb{R})=\mathbb{R}$, and the underlined $\mathbb{R}$ on the bottom right is the one corresponding to (\ref{nontrivial}). In both cases, the bottom of the left tower lies in degree $0$, and $U$ acts vertically.
\end{remark}

\begin{example}If the piercing sequence is $(+,-,+,-)$, we have
\begin{equation*}
\check{C}_*(Y,\spin)=\T^{\oplus 2}\oplus \T^{\oplus 2}\langle-1\rangle
\end{equation*}
where for each $k\geq 0$ in degrees $2k$ and $2k+1$ the chain complex looks like the Morse chain complex of the circle $S^1$ equipped with a Morse function with exactly two maxima and two minima. In particular
\begin{equation*}
\HMt_*(Y,\spin)=\T\oplus \T\langle-1\rangle\oplus\mathbb{Z}_{-1}
\end{equation*}
and the action of $\gamma$ is a surjective from the right tower onto the left tower, and zero on the bottom of the right tower and the summand $\mathbb{Z}_{-1}$. Correspondingly, we have for a local coefficient system $\Gamma_\eta$ with $[\eta]$ non-zero
\begin{equation*}
\HMred_*(Y,\spin;\Gamma_\eta)=\mathbb{R}^{\oplus 2}_{-1}.
\end{equation*}
Notice that in this situation the reduced Floer homology (with untwisted coefficients) $$\HMred_*(Y,\spin)=\mathbb{Z}_{-1}$$is also non-trivial.
\end{example}

\vspace{0.5cm}

\section{Trace formulas and families of Dirac operators}\label{traceformulas}

In this section we describe the main tools that we will employ to compute piercing sequences in our examples of interest. We will first recall the trace formulas for Dirac operators proved in \cite{LL3}, and then derive the trace formula that we will use to study the derivative of the eigenvalues.

\subsection{The geometric setup. } As in \cite[Section 2]{LL3}, we will fix once for all a spin structure on $Y$, which we think of as a lift
\begin{equation*}
\begin{tikzcd}
	& {SL(2;\mathbb{C})} \\
	{\pi_1(Y)} & {\mathbb{P}SL(2;\mathbb{C})}
	\arrow[from=1-2, to=2-2]
	\arrow[dashed, from=2-1, to=1-2]
	\arrow[hook, from=2-1, to=2-2]
\end{tikzcd}
\end{equation*}
after identifying $\mathrm{Isom}^+(\mathbb{H}^3)=\mathbb{P}SL(2;\mathbb{C})$.
\par
This lift determines for each closed geodesic $\gamma$ the spin holonomy $\hol(\tilde{\gamma})$, which is one of the two square roots of $\hol({\gamma})\in\mathbb{R}/2\pi\mathbb{Z}$. In \cite{LL3}, we showed how this data can be computed explicitly starting from a Dirichlet domain for $Y$ and a representative in $\mathbb{P}SL(2;\mathbb{C})$ for each closed geodesic. In the last step, we use the identification of closed geodesics with non-trivial conjugacy classes in $\pi_1(Y)$, cf. \cite{Mar}.
\par
Using this spin structure as a basepoint, each flat $B$ spin$^c$ connection corresponds to a twisting homomorphisms
\begin{equation*}
\varphi_B:\pi_1(Y)\rightarrow U(1).
\end{equation*}
The homomorphism factors through the abelianization $H_1(Y;\mathbb{Z})$ of $\pi_1(Y)$, and its restriction to the torsion subgroup determines the topological type of the corresponding torsion spin$^c$ structure. There is a corresponding notion of spin$^c$ holonomy, which at the practical level is obtained from the spin holonomy via the corresponding twisting homomorphism; to compute the latter, we only need to know the class $[\gamma]\in H_1(Y;\mathbb{Z})$ of the geodesic, which is also computed in \cite{LL3} taking as input a Dirichlet domain for $Y$.

\vspace{0.5cm}

\subsection{The trace formulas.} For a fixed torsion spin$^c$ structure $\spin$, we will consider the family of Dirac operators
\begin{equation*}
D_{B_{\tau}}:\Gamma(S)\rightarrow \Gamma(S)
\end{equation*}
for $[B_\tau]\in\mathbb{T}$ varying in the torus of flat spin$^c$ connection for $(Y,\spin)$. We parametrize the latter (after choosing a basepoint $B_0$) by $\tau\in[0,1]$. Denote by
\begin{equation}\label{eigenvaluesdir}
\dots\leq s_{-2}(\tau)\leq s_{-1}(\tau)\leq s_0(\tau)\leq s_1(\tau)\leq s_2(\tau)\leq \dots
\end{equation}
the eigenvalues of $D_{B_{\tau}}$ with repetitions, and by $\varphi_\tau$ the corresponding twisting homomorphism
\begin{remark}
By the Weyl Law \cite{Roe}, we have
\begin{equation}\label{weyl}
\#\{j\text{ such that }|s_j(\tau)|\in[T,T+1]\}\sim C\cdot T^2
\end{equation}
for some constant $C>0$ as $T\rightarrow \infty$.    
\end{remark}
Consider test functions $H,K$ which are compactly supported and even and odd respectively. Denote by $\widehat{H},\widehat{K}$ their Fourier transforms, for which we use the convention
\begin{equation}\label{fourierconvention}
\widehat{H}(t)=\int_{\mathbb{R}}H(x)e^{-itx}dx, 
\end{equation}
and similarly for $K$.
\begin{remark}
Throughout the paper, we will use $\tau$ to parametrize $\mathbb{T}$ and $t$ for the variable in frequency space.
\end{remark}
The functions $\widehat{H}$ and $\widehat{K}$ are real and purely imaginary respectively. As in \cite{LL1}, will assume that $H,K$ are regular enough, meaning that
\begin{equation}
\int_{\mathbb{R}}\left(|\widehat{H}(t)|^2+|\widehat{H}'(t)|^2\right) \left(\sqrt{1+t^2})\right)^{2\delta}<\infty
\end{equation}
for some $\delta>5/2$ and similarly for $K$; this is of course true if $H,K$ are smooth, but for our purposes it will be convenient to use less regular functions.
\par
Then the identities
\begin{equation}\label{tracediraceven}
\underbrace{\frac{1}{2}\sum \widehat{H}(s_j(\tau))}_{\text{spectral side}}= \underbrace{\frac{\mathrm{vol}(Y)}{2\pi}\left(\frac{1}{4}H(0)-H''(0)\right)+\sum \ell(\gamma_0)\frac{\cos(\hol(\tilde{\gamma}))\cdot\cos(\varphi_\tau(\gamma))}{|1-e^{\mathbb{C}\ell(\gamma)}||1-e^{-\mathbb{C}\ell(\gamma)}|}  H(\ell(\gamma))}_{\text{geometric side}}
\end{equation}
and
\begin{equation}\label{tracediracodd}
\underbrace{-\frac{i}{2}\sum \widehat{K}(s_j(\tau))}_{\text{spectral side}}=\underbrace{\sum \ell(\gamma_0)\frac{\sin(\hol(\tilde{\gamma}))\cdot \cos(\varphi_\tau(\gamma))}{|1-e^{\mathbb{C}\ell(\gamma)}||1-e^{-\mathbb{C}\ell(\gamma)}|}  K(\ell(\gamma))}_{\text{geometric side}}
\end{equation}
hold for each $\tau\in[0,1]$, see \cite[Section 3]{LL3}; again the sum on the spectral side is taken over all eigenvalues, while the sum over the geometric side is taken over all closed geodesics.
\vspace{0.3cm}
\subsection{The derivative of the trace formula.} For the purposes of proving transversality we will need to show that the derivative of the smallest eigenvalue is non-vanishing in suitable intervals. In order to obtain information about such derivative, we will differentiate the \text{odd} trace formula (\ref{tracediracodd}) as follows. Having $H^1(Y)=\mathbb{Z}$ implies that
\begin{equation*}
\varphi_\tau(\gamma)=e^{2\pi i\cdot \tau[\gamma]}\cdot\varphi_0(\gamma)
\end{equation*}
where $[\gamma]$ is the class in $H_1(Y)/\mathrm{tors}=\mathbb{Z}$.  Differentiating formally in the variable $\tau$ the odd trace formula (\ref{tracediracodd}), we obtain by the chain rule the expression
\begin{equation}\label{derivativetraceodd}
\underbrace{-\frac{i}{2}\sum \widehat{K}'(s_j(\tau))s'_j(\tau)}_{\text{spectral side}}=\underbrace{2\pi  \sum \ell(\gamma_0)\frac{\sin(\hol(\tilde{\gamma}))}{|1-e^{\mathbb{C}\ell(\gamma)}||1-e^{-\mathbb{C}\ell(\gamma)}|}\sin(\varphi_\tau(\gamma))\cdot [\gamma]\cdot K(\ell(\gamma))}_{\text{geometric side}}.
\end{equation}
\begin{remark}
Of course, this can be done for the family of even trace formulas as well; for our purposes of estimating the derivative of the small eigenvalues we consider the odd case because for $H$ even $\widehat{H}'(0)=0$ so that the contribution of the derivative of small eigenvalues is negligible.
\end{remark}
There are two technical points to discuss in order to make the formula rigorous. As a starting point for both of them, notice that the family of Dirac operators $\{D_{B_\tau}\}$ can be written (up to gauge equivalence) in the form
\begin{equation}\label{diracfamily}
D_{B_\tau}=D_{B_0}+2\pi i\tau\cdot \rho(\alpha)\text{ for }\tau\in[0,1]
\end{equation}
where $\rho$ is Clifford multiplication, and $\alpha$ is any closed $1$-form representing the generator of $H^1(Y;\mathbb{Z})$.
\\
\par
First of all, Rellich's theorem (see for example \cite[Ch. 14.9]{Bus}) implies that the $1$-parameter family of eigenvalues can be \textit{locally} written uniquely as the graph of analytic (and in particular continuously differentiable) functions $\tilde{s}_j$. Notice though that in general for $\tau\neq \tau_0$ the parametrization is different from the one given by the natural ordering (\ref{eigenvaluesdir}); the local example to keep in mind is
\begin{equation*}
s_0(\tau)=-|\tau|\text{ and }s_1(\tau)=|\tau|,
\end{equation*}
for which the analytic parametrizations are
\begin{equation*}
\tilde{s}_0(\tau)=\tau\text{ and }\tilde{s}_1(\tau)=-\tau.
\end{equation*}
To avoid cluttering the notation we will denote both parametrization by $\{s_i(\tau)\}$; this will not cause confusion at any point in the paper.
\par
While the discussion in \cite{Bus} focuses on the case of the family of Laplacians on hyperbolic surfaces (parametrized by Teichm\"uller space), the proof readily adapts to our setup because our family of Dirac operators (locally) extends to the family of operators
\begin{equation*}
T_z=D_{B_0}+iz\cdot 2\pi\rho(\alpha) \text{ for }z\in\mathbb{C}
\end{equation*}
which is \textit{holomorphic} in the sense that for fixed spinors $\Phi,\Psi \in \Gamma(S)$, the map
\begin{equation*}
z\mapsto\langle T_z\Phi,\Psi\rangle
\end{equation*}
is holomorphic in $z$; in fact our situation is much simpler than the one arising in \cite{Bus} as we are dealing with a linear family.
\begin{remark}
Notice also that
Rellich's theorem is false for higher dimensional families, already in finite dimensions \cite[Ch. 14.9]{Bus}.
\end{remark}
Second, basic spectral pertubation theory (see for example \cite[II.3]{Kato}) for a family of self-adjoint operators on a Hilbert space of the form $A+\tau B$ with $B$ bounded provides a bound on the derivative of the eigenvalues in terms of the operator norm of $B$. In our case the underlying Hilbert space is $L^2(S)$ and because Clifford multiplication
\begin{equation*}
\rho: T^*Y\rightarrow \mathfrak{su}(S)
\end{equation*}
is an isometry of Euclidean bundles, where on the latter we consider the inner product $\mathrm{tr}(a^*b)/2$, the operator norm of $2\pi i \cdot\rho(\alpha)$ is simply $2\pi\|\alpha\|_{C^0}$. Setting
\begin{equation}
C_Y=\inf_{\alpha}\|\alpha\|_{C^0(Y)}
\end{equation}
where the infimum is taken over the closed $1$-forms $\alpha$ with $[\alpha]\in H^1(Y;\mathbb{Z})$ a generator, we therefore have the following a priori bound.
\begin{lemma}\label{eigenbound}
For all eigenvalues $s_j(\tau)$ (in the analytic parametrization) the derivative satisfies the bound $|s'_j(\tau)|\leq 2\pi C_Y$.
\end{lemma}
From this, using the Weyl law (\ref{weyl}) we readily conclude that the formula (\ref{derivativetraceodd}) holds for $K$ smooth and compactly supported. The same approximation argument as in \cite[Appendix $C$]{LL1} proves that the following holds.
\begin{prop}\label{derivativetrace}
Consider a compactly supported odd test function $K$ which is not necessarily smooth but satisfies
\begin{equation}\label{norm}
\int_{\mathbb{R}}\left(|\widehat{K}(t)|^2+|\widehat{K}'(t)|^2+|\widehat{K}''(t)|^2\right) \left(\sqrt{1+t^2})\right)^{2\delta}<\infty
\end{equation}
for some $\delta>5/2$. Then for each $\tau\in[0,1]$ the trace formula (\ref{derivativetraceodd}) holds.
\end{prop}

The proof is indeed the same as in \cite{LL1} and follows from the boundedness under this norm of the linear functional
\begin{equation*}
K\mapsto \sum \widehat{K}'(s_j(\tau))s'_j(\tau)
\end{equation*}
which, using Lemma (\ref{eigenbound}), is proved by showing that the function
\begin{equation*}
K\mapsto \sum |\widehat{K}'(s_j(\tau))|
\end{equation*}
is bounded in terms of the norm (\ref{norm}).

\begin{remark}
For the purposes of the proof of our main results, we will need explicit upper bounds on $C_Y$ for our examples; this is done is Appendix \ref{C0bounds}.
\end{remark}

\vspace{0.5cm}

\section{Identifying intervals with small eigenvalues}\label{smalleigint}

\begin{figure}
\centering
  \includegraphics[width=.6\linewidth]{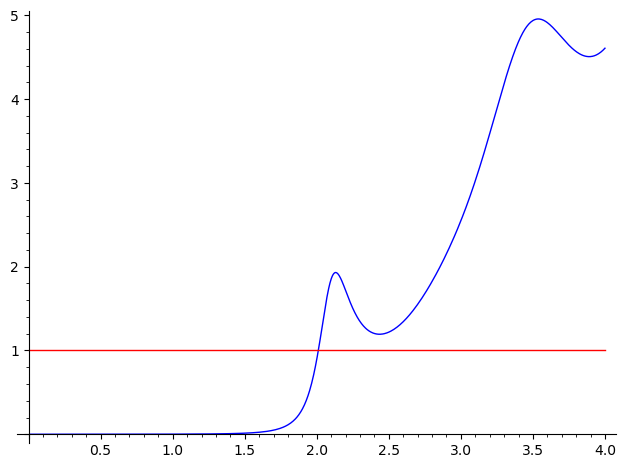}
  \caption{The output of the techniques developed in \cite{LL1} for Census 357 at cutoff $R=6.5$: if $\lambda$ is an eigenvalue of the Hodge Laplacian on coexact $1$-forms, then the blue plot is $\geq1$ at $\sqrt{\lambda}$. In particular, $\lambda_1^*>2$.}
    \label{357coexact}
\end{figure}

We begin the discussion of our strategy to show that for a given hyperbolic three-manifold equipped with torsion spin$^c$ structure $(Y,\spin)$, Condition $2$ holds. Notice that this is an independent task from checking that Condition $1$ holds. As a concrete example, we will focus on the unique self-conjugate spin$^c$ structure on $\#357$ (which  has $H_{1}(Y;\mathbb{Z})_{\mathrm{tors}}=\mathbb{Z}/7\mathbb{Z}$), which is readily checked to be spectrally large, see Figure \ref{357coexact}. 
\par
The determination of the intervals in $\mathbb{T}$ for which the Dirac operator $D_B$ admits a small eigenvalue involves two steps: first, determining intervals potentially containing small eigenvalues, and then certifying that the corresponding operators actually admit them. We discuss the two parts separately.

\vspace{0.3cm}

\subsection{Potential intervals with small eigenvalues.}\label{formalbooker} The method introduced in \cite{LL1}, applied to a specific value of $\tau\in\mathbb{T}$, allows one to determine explicitly a function $J_s(\tau)\geq 0,$ for $s\geq 0,$ with the property that
\begin{equation}\label{bookerbound}
J_s(\tau)\geq \text{multiplicity of $s$ as absolute value of eigenvalues of $D_{B_\tau}$}.
\end{equation}
Notice that the outcome depends on some choices such as a cutoff $R> 0$ for the length spectrum and a finite dimensional vector space $\mathcal{F}$ of test functions, which in our implementation has as basis suitable shifts of the convolution square of the indicator function of an interval.
\par
As a special case of (\ref{bookerbound}) we see that if $J_s(\tau) < 1$, then $s$ is not an eigenvalue of $D_{B_\tau}$.
\par
In particular, for the locus $\mathsf{K}$ of operators with kernel, the containment
\begin{equation}\label{smalleiglocus}
\mathsf{K}\subset\{\tau\text{ such that } J_0(\tau)\geq 1\}
\end{equation}
holds. We expect the latter to generically be a union of closed disjoint intervals in $\mathbb{T}$, and its determination can be performed very efficiently by exploiting the specific structure of the $1$-parameter family of trace formulas (\ref{tracediraceven}), as we now describe.
\par
Let us first recall that for a fixed $\tau\in\mathbb{T}$, one obtains an $n\times n$ positive definite matrix $A(\tau)$ (where $n=\dim\mathcal{F}$) whose entries are obtained by evaluating the geometric side trace formula (\ref{tracediraceven}) at suitable functions which are convolutions of elements in $\mathcal{F}$. We then have
\begin{equation*}
J_s(\tau)=\frac{1}{\langle\ v_s, A(\tau)^{-1}v_s\rangle}
\end{equation*}
where $v_s$ is an $n$-column vector of functions of the spectral variable $s$ explicitly determined in terms of $\mathcal{F}$; in our specific implementation, these are elementary functions involving sines and cosines.
\\
\par
With the goal of letting $\tau$ vary in $[0,1]$, one can then compute for any even compactly supported test function $H$ a `formal' geometric side
\begin{equation}\label{formalgeometricside}
\frac{\mathrm{vol}(Y)}{2\pi}\left(\frac{1}{4}H(0)-H''(0)\right)+\sum \ell(\gamma_0)\frac{\cos(\hol(\tilde{\gamma}))}{|1-e^{\mathbb{C}\ell(\gamma)}||1-e^{-\mathbb{C}\ell(\gamma)}|}  H(\ell(\gamma))\cdot [\gamma]\in \mathbb{C}[H_1(Y;\mathbb{Z})]
\end{equation}
with values in the group ring $\mathbb{C}[H_1(Y;\mathbb{Z})]$; this sum is supported on the finitely many closed geodesics of length at most $R$ and hence is well-defined. The value for the geometric side at a specific parameter $\tau\in\mathbb{T}$ is then obtained by applying the algebra homomorphism 
\begin{equation*}
\mathbb{C}[H_1(Y;\mathbb{Z})]\rightarrow \mathbb{C}
\end{equation*}
corresponding to the twisting homomorphism $\varphi_\tau$ (which factors through $H_1(Y;\mathbb{Z})$).
\par
Going back to the eigenvalue bounds, combining this approach with the methods introduced in \cite{LL1} allows us to determine a `formal' $n\times n$ matrix $A$ with entries in $\mathbb{C}[H_1(Y;\mathbb{Z})]$, such that each given $A(\tau)$ is gotten by evaluating $A$ via the homomorphism $\varphi_\tau$. An important consequence of this description is that the dependence of the matrix $A(\tau)$ in $\tau$ is through elementary functions involving sines and cosines, so $J_s(\tau)$ has the same dependence as well.  This allows us to compute the locus (\ref{smalleiglocus}) accurately and quickly.  

\begin{figure}
\centering
  \includegraphics[width=.6\linewidth]{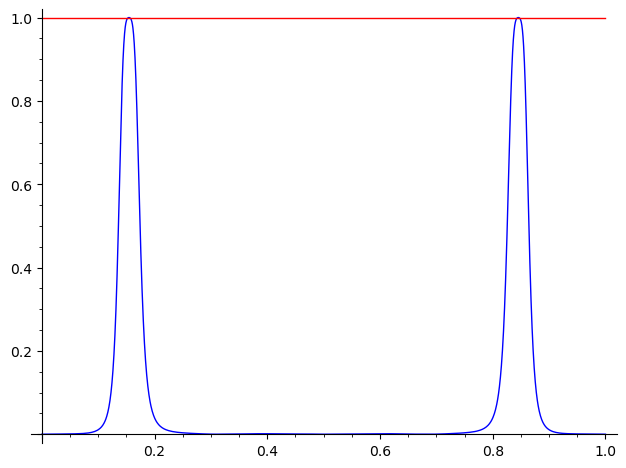}
  \label{fig:sub1}
  \caption{The plot of $J_0(\tau)$ for the unique self-conjugate spin$^c$ structure on Census 357, computed at cutoff $R=7$.}
  \label{J0357spin}
\end{figure}

For example, for the unique self-conjugate spin$^c$ structure on the census manifold \#357, the plot of the function $J_0(\tau)$ for $\tau\in[0,1]$ based on the computation of the length spectrum up to cutoff $R=7$, is shown in Figure \ref{J0357spin}. The intervals potentially containing points of $\mathsf{K}$ are 
\begin{equation*}
[0.1537,0.1556]\cup[0.8444, 0.8463]
\end{equation*}
respectively. As expected, these are symmetric under the conjugation symmetry $\tau\mapsto 1-\tau$ for self-conjugate spin$^c$ structures.

\vspace{0.3cm}
\subsection{Certifying a unique small eigenvalue.} 
\begin{figure}
\centering
  \includegraphics[width=.6\linewidth]{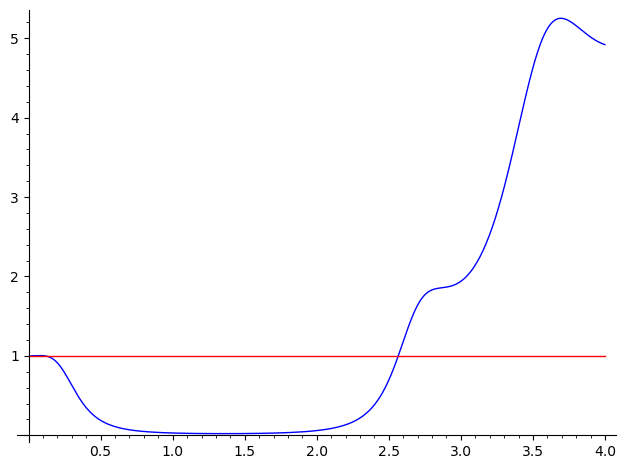}
  \label{fig:sub1}
  \caption{The plot of $J_s(\tau_0)$ at $\tau_0=0.15465$ for the unique self-conjugate spin$^c$ structure on Census 357, computed at cutoff $R=7$.}
  \label{J0357spinmidpoint}
\end{figure}

We now show that for our examples, in each of the intervals in $\mathbb{T}$ determined above the Dirac operator admits \textit{exactly} one small eigenvalue; of course, by symmetry we can focus on the first of the two $[0.1537,0.1556]$. One can plot the functions $J_s(\tau)$ for a fixed value of $\tau$; the plots for the midpoint of the interval $\tau_0=0.15465$ obtained from discussion above is shown in Figure \ref{J0357spinmidpoint}.

\vspace{0.3cm}

Focusing first at time $\tau_0$, we have that $J_s(\tau_0)\geq1$ for $s$ small, with the graph crossing the value $1$ at some point in the interval $[0.0736,0.0737]$; furthermore the next spectral parameter is larger that $2.5$. This is clearly suggestive of the fact that there is exactly one small eigenvalue. To certify that this is actually the case, we will show separately that $D_{B_{\tau_0}}$ admits at most one eigenvalue with absolute value $\leq 0.0737$, and then show that there is at least one. These are both achieved by applying the trace formula (\ref{tracediraceven}) for suitable even test functions as follows.
\\
\par
\textit{Step 1: Showing that $D_{B_{\tau_0}}$ has at most one small eigenvalue.} This can be achieved by looking at test even test functions $H$ for which $\widehat{H}(t)\geq 0$ and $\widehat{H}(0)$ is somewhat large. In practice, a good choice for this example is
\begin{equation*}
H_6(x)=\left(\frac{1}{2}\mathbf{1}_{[-1,1]}\right)^{\ast 6}(x)
\end{equation*}
which is supported in $[-6,6]$ and has $\widehat{H_6}(x)=\sinc^6(x)$, where
\begin{equation*}
\sinc(x)=\begin{cases}
\frac{\sin(x)}{x} \text{ if }x\neq0\\
0 \text{ otherwise.}
\end{cases}
\end{equation*}
For $L\leq 0.5$, denoting by $\Gamma_{H_6}(\tau_0)$ the geometric side of (\ref{tracediraceven}) evaluated with $H_6$, we have that
\begin{equation*}
\Gamma_{H_6}(\tau_0)=\frac{1}{2}\sum_{s_j(\tau_0)}\widehat{H_6}(s_j(\tau_0))\geq \frac{1}{2}\sum_{|s_j(\tau_0)|\leq L}\widehat{H_6}(s_j(\tau_0))
\end{equation*}
because $\widehat{H_6}\geq0$, so that we conclude
\begin{equation*}
\#\{\text{spectral parameters with }|s_j(\tau_0)|\leq L\}\leq \left\lfloor\frac{2\cdot\Gamma_{H_6}(\tau_0)}{\min_{x\in[0,L]}\sinc^6(x)}\right\rfloor.
\end{equation*}
In our situation, we conclude that the Dirac operator $D_{B_{\tau_0}}$ admits at most one eigenvalue with absolute value $\leq0.0737$.
\\
\par
\textit{Step 2: Showing that $D_{B_{\tau_0}}$ has at least one small eigenvalue.} This can achieved, after choosing an even test functions $H$ for which $\widehat{H}(t)\geq 0$ and $\widehat{H}(0)$ is somewhat large, by considering the trace formula evaluated at the function
\begin{equation*}
G(x)=H(x)-\frac{H''(x)}{L^2},
\end{equation*}
which has Fourier transform
\begin{equation*}
\widehat{G}(t)=\widehat{H}(t)\cdot \left(1-\frac{t^2}{L^2}\right).
\end{equation*}
In particular, $\widehat{G}(t)\leq 0$ if $|t|\geq L$, so that if there were no eigenvalues in $[-L,L]$ we would have $\Gamma_G(\tau_0)\leq 0$, where $\Gamma_G(\tau_0)$ denotes the geometric side of (\ref{tracediraceven}) evaluated using the function $G$.
\par
In our particular situation, the discussion applied to the test function $H=H_6$ shows that there is at least one eigenvalue with absolute value $\leq0.0737$. 
\\
\par
Finally, from our computation at the midpoint $\tau_0$, we conclude that for $\tau\in[0.1537,0.1556]$ the Dirac operator $D_{B_{\tau}}$ admits exactly one small eigenvalue by simply looking at the plots of $J_s(\tau)$ (which change very little from $J_s(\tau_0)$), and the continuity of the eigenvalues $s_i(\tau)$; alternatively, we can repeat the arguments of Step $1$ and $2$ above at other the points in the interval. In what follows, we will denote the small eigenvalue by $s_0(\tau)$.

\vspace{0.5cm}

\section{Certifying single transverse crossings}\label{singlecrossing}

In this section we conclude the proof of Theorem \ref{thm1} for the unique self-conjugate spin$^c$ structure on \#357 by showing that the intervals identified in the previous section contain a single, transverse crossing, so that the piercing sequence is $(+,-)$.
\par
Notice that in the previous section we have only worked with even test functions, and therefore we have gathered no information about the sign of the small eigenvalue $s_0(\tau)$. The key point of this section is then to employ the trace formulas (\ref{tracediracodd}) and (\ref{derivativetraceodd}) involving odd test functions in order to study the latter; in particular, we will treat the terms on the spectral side involving eigenvalues $\{s_i\}_{i\neq 0}$ as an error term to be suitably bounded.
\par
We will recall how to explicitly bound the number of eigenvalues in a given interval in the first subsection. Then, in the following two subsections we prove the existence of a crossing and its uniqueness and transversality, respectively.

\vspace{0.3cm}
\subsection{Upper bounds on spectral densities. }\label{upperbounddensity}For $T\geq0$, we want to provide concrete bounds on the number of eigenvalues with absolute value in a given interval $[T,T+1]$; this is the content of explicit local Weyl laws as obtained in \cite{LL3} using the even trace formula (\ref{tracediraceven}). We will use the test function
\begin{align*}
H_{6,\nu}&=\left(\frac{1}{2}\mathbf{1}_{[-1,1]}\right)^{*6}\cdot(e^{i\nu x}+e^{-i\nu x})=\\
&=2\cdot \left(\frac{1}{2}\mathbf{1}_{[-1,1]}\right)^{*6}\cdot\cos(\nu x).
\end{align*}
for which
\begin{align*}
-H''_{6,\nu}(0)&=\frac{11}{20}\nu^2+\frac{1}{4}\\
\widehat{H_{6,\nu}}(t)&=\sinc^6(t+\nu)+\sinc^6(t-\nu).
\end{align*}
Using that
\begin{itemize}
\item $\widehat{H_{6,\nu}}(t)\geq 0$ for all $t$;
\item $\sinc^6(t)\geq 0.777$ for $|t|\leq 1/2$,
\end{itemize}
we obtain as in the previous section that the inequality
\begin{equation}\label{spectralparameterbound}
\#\left\{\text{spectral parameters with }|s_j(\tau)|\in[\nu-1/2,\nu+1/2]\right\}\leq \left\lfloor\frac{2}{0.777}\cdot\Gamma_{H_{6,\nu}}(\tau)\right\rfloor
\end{equation}
holds for $\nu\geq 1/2$. This bound is readily evaluated as $\tau$ varies using the method described in Section \ref{formalbooker}.
\begin{remark}\label{mainasymptotic}
When $\nu$ is large, the main contribution to the geometric side of the Selberg trace formula arises from the quadratic term in $\nu$ in the second derivative of $H''_{6,\nu}$.
\end{remark}

\vspace{0.3cm}

\subsection{Proving the existence of a crossing. }
\begin{figure}
\centering
  \includegraphics[width=.5\linewidth]{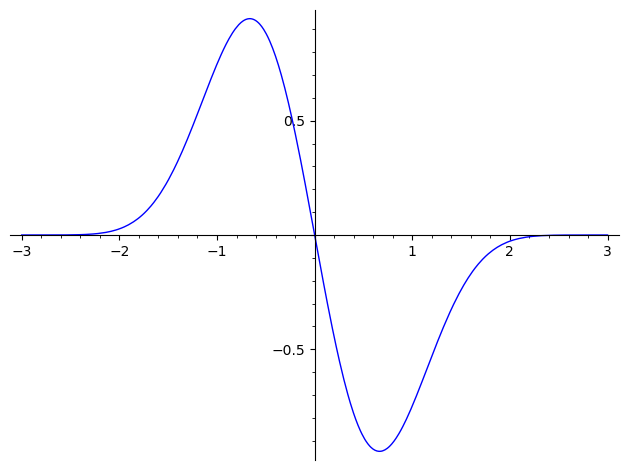}
  \label{fig:sub1}
  \caption{The plot of the first derivative of $\sinc^7$. This function is useful to detect the sign of `medium' sized eigenvalues.}
  \label{sinc7firstder}
\end{figure}

A major complication when dealing with the trace formula for a odd test function $K$ is that $\widehat{K}(0)=0$, so that it is challenging to understand the contribution of small eigenvalues. The strategy to certify the existence of a crossing in the interval is then to look at nearby parameters \textit{without} small eigenvalues. More specifically, we will consider the odd test function
\begin{equation}\label{xconv7}
K_7(x)=x\cdot \left(\frac{1}{2}\mathbf{1}_{[-1,1]}(x)\right)^{*7}
\end{equation}
which, according to our convention (\ref{fourierconvention}), has
\begin{equation*}
-i\widehat{K_7}(t)=(\sinc^7)'(t).
\end{equation*}
Explicitly, we have for $t\neq 0$ that the first derivative is
\begin{equation*}
(\sinc^7)'(t)=\frac{7\cdot\sin^6(t)\cdot(t\cdot \cos(t)-\sin(t))}{t^8};
\end{equation*}
to help the reader's intuition, the plot of this function is shown in Figure \ref{sinc7firstder}.

\vspace{0.3cm}

\begin{figure}
\centering
  \includegraphics[width=.6\linewidth]{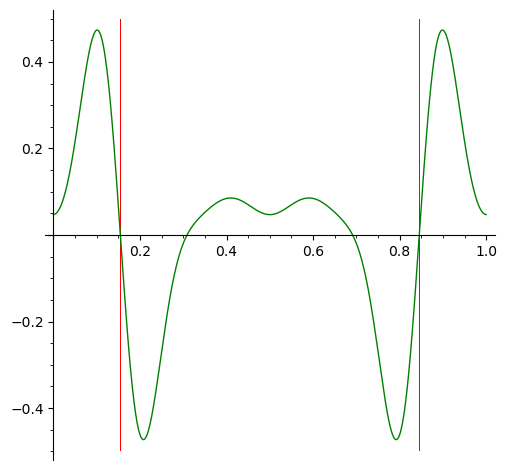}
  \label{fig:sub1}
  \caption{In green, the plot of $\Gamma_{K_7}(\tau)$ for unique self-conjugate spin$^c$ structure on Census \#357. The vertical red lines are really bars denoting the interval $[0.1537,0.1556]$ and its mirror.}
  \label{357spinoddtraceformulapic}
\end{figure}

In order to get a heuristic understanding of the sign of the eigenvalues, we first consider the geometric side $\Gamma_{K_7}(\tau)$ of the odd trace formula (\ref{tracediracodd}) at parameter $\tau\in\mathbb{T}$; the plot can be found in Figure \ref{357spinoddtraceformulapic}. This is readily evaluated by adapting with the method described for even functions in Subsection \ref{formalbooker} to the odd case. Focusing on the first interval, the fact that $\Gamma_{K_7}(\tau)$ goes from very positive to very negative is quite suggestive of the fact that there is a crossing in the interval $[0.1537,0.1556]$.

\vspace{0.3cm}

\begin{figure}
\centering
\begin{subfigure}{.45\linewidth}
    
  \includegraphics[width=\linewidth]{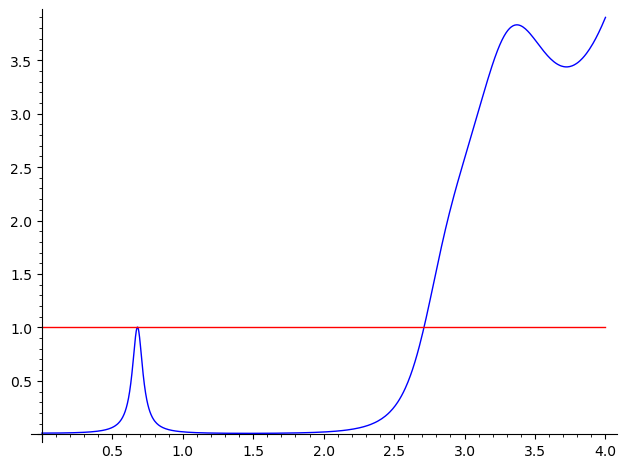}
  
\end{subfigure}
\begin{subfigure}{.45\linewidth}
  \includegraphics[width=\linewidth]{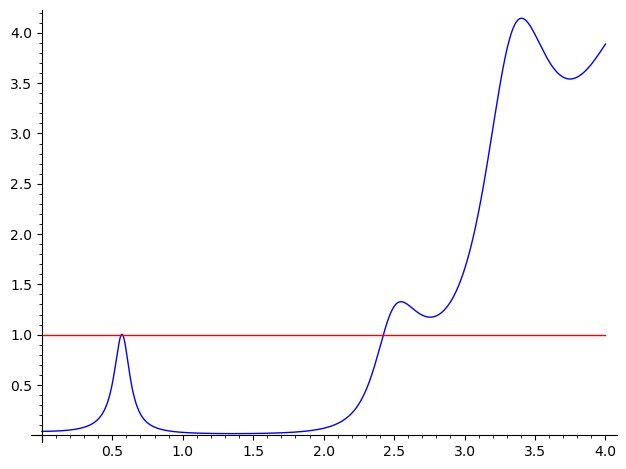}

  \end{subfigure}
    \caption{The plots of $J_s(0.1)$ and $J_s(0.2)$ for the unique self-conjugate spin$^c$ structure on Census 357, computed at cutoff $R=7$.}

      \label{357spinbooker0102}
\end{figure}

In order to make this rigorous, we consider the specific values $\tau=0.1$ and $\tau=0.2$, the corresponding plots of $J_s(\tau)$ can be found in Figure \ref{357spinbooker0102}. Looking at the whole family of functions $J_s(\tau)$ for $\tau\in[0.1,0.2]$, we see that the `peaks' above $1$ of the functions $J_0(0.1)$ and $J_0(0.2)$ fit in a family that for $\tau$ in the interval $[0.1537,0.1556]$ is the peak near $0$ as in Figure \ref{J0357spinmidpoint}. In particular, they correspond to the eigenvalue $s_0(\tau)$ for $\tau\in[0.1,0.2]$.
\par
Furthermore, we evaluate that
\begin{equation*}
\Gamma_{K_7}(0.1)\approx0.4735 \text{ and }\Gamma_{K_7}(0.2)\approx -0.4616.
\end{equation*}
Because of the odd trace formula (\ref{tracediracodd}), we have
\begin{equation*}
\Gamma_{K_7}(\tau)=\frac{1}{2}\sum_{i\in\mathbb{Z}}(\sinc^7)'(s_i(\tau)),
\end{equation*}
so that this proves $s_0(0.1)<0$ and $s_0(0.2)>0$ provided we can suitably bound the contribution of the other eigenvalues. By continuity, we then conclude $s_0(\tau)$ attains the value $0$ for some $\tau$ in our interval $[0.1537,0.1556]$.

\vspace{0.3cm}

Regarding $\tau=0.1$, using that the next eigenvalue has absolute value at least $2.7$, we have using the spectral density bounds in (\ref{spectralparameterbound}) that
\begin{align*}
\left|\sum_{i\neq0}(\sinc^7)'(s_i(0.1))\right|&\leq\sum_{k\geq 0}\max_{t\in[2.7+k,2.7+k+1]}|(\sinc^7)'(t)|\cdot \#\left\{|s_i(0.1)|\in[2.7+k,2.7+k+1]\right\}\\
&\leq 0.0017
\end{align*}
so that the conclusion follows.

\begin{remark}\label{numericalapprox}
At the practical level, for simplicity we bound the infinite sum by evaluating its first forty terms. This does not introduce a significant error because we have$$(\sinc^7)'(t)=O(t^{-7}),$$while the spectral density bounds are $O(t^2)$. This tail can be bounded explicitly using bounding (using $|\cos(x)|\leq1)$ the sum of $H_{6,\nu}$ over the closed geodesics with
\begin{equation*}
S_Y=\sum \frac{\ell(\gamma_0)}{|1-e^{\mathbb{C}\ell(\gamma)}||1-e^{-\mathbb{C}\ell(\gamma)}|}H_{6,0}(\ell(\gamma))
\end{equation*}
which is readily computable, see also Remark \ref{mainasymptotic}.
\end{remark}
At $\tau=0.2$, we have that the next eigenvalue has absolute value at least $2.4$ and obtain correspondingly that
\begin{equation*}
\left|\sum_{i\neq0}(\sinc^7)'(s_i(0.2))\right|\leq 0.0104,
\end{equation*}
and conclude again.

\vspace{0.3cm}

\subsection{Proving uniqueness and transversality.}
\begin{figure}
\centering
  \includegraphics[width=.5\linewidth]{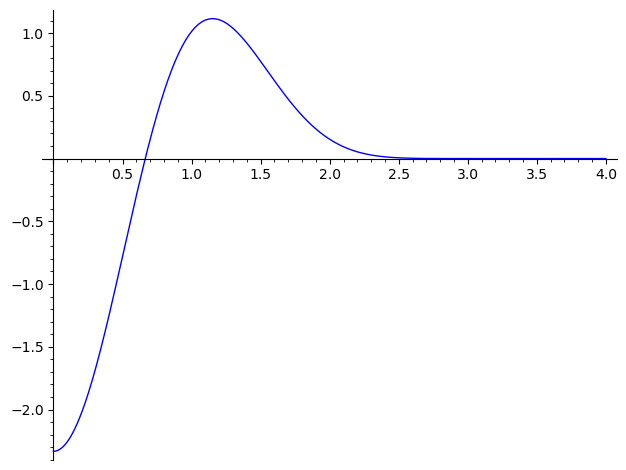}
  \label{fig:sub1}
  \caption{The plot of the second derivative of $\sinc^7$ for $t\geq0$. This function is very small for $t \gtrsim 2.5$.}
  \label{sinc7secondder}
\end{figure}

We are left to show that for the small eigenvalue
\begin{equation*}
s_0'(\tau)\neq 0\text{ for }\tau\in [0.1537,0.1556].
\end{equation*}
In order to do this, we will study the derivative of the trace formula (\ref{derivativetraceodd}) for the odd test function $K_7$ in (\ref{xconv7}). We have  

\begin{equation*}
-i\widehat{K_7}'(t)=(\sinc^7)''(t)
\end{equation*}
where explicitly the second derivative is
\begin{equation*}
(\sinc^7)''(t)=-\frac{7\sin^5(t)\cdot\left(\sin(t)\cdot((t^2-2)\sin(t)+2t\cdot\cos(t))-6(\sin(t)-t\cos(t))^2\right)}{t^9}
\end{equation*}
for $t\neq0$, with value $-7/3$ at zero; see Figure \ref{sinc7secondder} for the plot.

\vspace{0.3cm}

For our purposes, we will rewrite the trace formula (\ref{derivativetraceodd}) as
\begin{equation*}
-i\widehat{K_7}'(s_0(\tau))s'_0(\tau)=2\cdot\widetilde{\Gamma}_{K_7}(\tau)+i \sum_{j\neq0} \widehat{K_7}'(s_j(\tau))s'_j(\tau)
\end{equation*}
where $\widetilde{\Gamma}_{K_7}(t)$ denotes the geometric side and we think of the infinite sum on the right hand side as an error term; in particular our goal will be to show that the inequality
\begin{equation}\label{desiredbound}
\left|\sum_{j\neq 0}\widehat{K_7}'(s_j(\tau))s'_j(\tau)\right|< 2|\widetilde{\Gamma}_{K_7}(\tau)|
\end{equation}
holds for $\tau\in [0.1537,0.1556]$, which implies $s_0'(\tau)\neq0$ in the interval.
\\
\par
Looking for example at the midpoint of the interval $\tau_0=0.15465$, we have
\begin{equation*}
\widetilde{\Gamma}_{K_7}(\tau_0)\approx-14.4905;
\end{equation*}
in fact, we have
\begin{equation*}
\widetilde{\Gamma}_{K_7}(\tau)\in[-14.4906,-14.4832]\text{ for all }\tau\in [0.1537,0.1556].
\end{equation*}
Using that the next eigenvalue of $D_{B_\tau}$ has absolute value at least $2.5$ throughout the interval (which is readily checked by evaluating $J_{2.5}(\tau)$), we see as in the previous subsection that
\begin{align*}
\left|\sum_{i\neq0}\widehat{K_7}'(s_i(\tau))\right|&\leq\sum_{k\geq 0}\max_{t\in[2.5+k,2.5+k+1]}|(\sinc^7)''(t)|\cdot \#\{|s_i(\tau)|\in[2.5+k,2.5+k+1]\}\\
&\leq 0.0470
\end{align*}
for all parameters $\tau$ in the interval. Again, we estimate the sum as in Remark \ref{numericalapprox}. As explained in Appendix \ref{C0bounds}, we have $C_Y\leq 3.5151$ for \#357,  so that by Lemma \ref{eigenbound} the we have the bound
\begin{equation*}
|s_i'(\tau)|\leq 2\pi C_Y\leq 22.0861\text{ for all }\tau\in[0,1],
\end{equation*}
so that inequality (\ref{desiredbound}) holds as
\begin{equation*}
0.0470\times 22.0861\leq 1.0381.
\end{equation*}
This concludes the proof of Theorem \ref{thm1} for the self-conjugate spin$^c$ structure on \#357.

\vspace{0.5cm}

\section{More examples}\label{moreexamples}

We have discussed the specific details of the proof of Theorem \ref{thm1} in the case of the self-conjugate spin$^c$ structure on \#357. We chose this example because in this situation the implementation of the methods are especially simple. We present in this section other spectrally large examples which are representative of the complications that arise when computing the piercing sequence in other examples.
\vspace{0.3cm}
\subsection{Non self-conjugate spin$^c$ structures on \#357. }
\begin{figure}
\centering
\begin{subfigure}{.45\linewidth}
  \includegraphics[width=\linewidth]{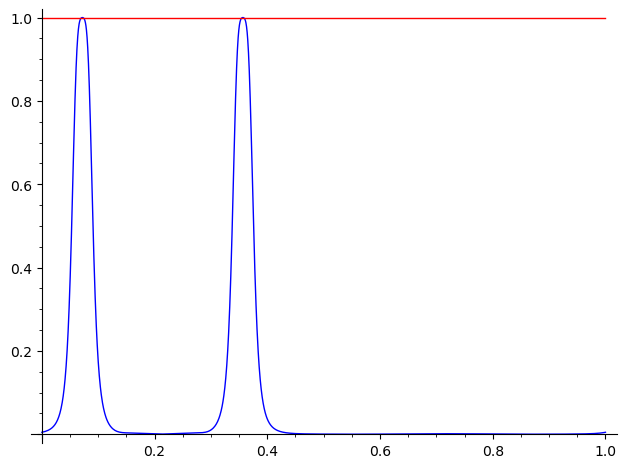}

\end{subfigure}
  \begin{subfigure}{.45\linewidth}
  \includegraphics[width=\linewidth]{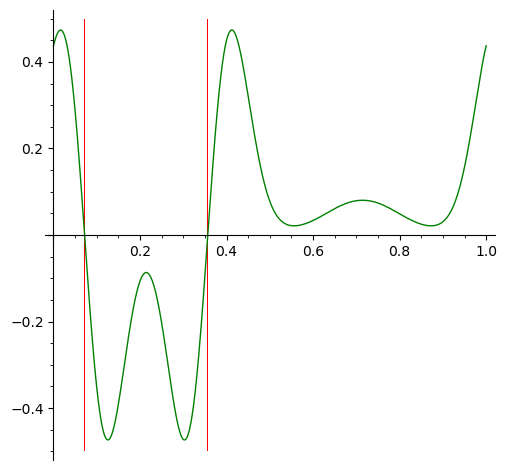}
 
  \end{subfigure}
   \caption{The plots of $J_0(\tau)$ and $\Gamma_{K_7}(\tau)$ for the pair of self conjugate spin$^c$ structures under consideration on \#357, computed at cutoff $R=7$.}
    \label{357nonspin}
\end{figure}

Other than the self-conjugate one, \#357 also admits a pair of conjugate spin$^c$ structures admitting non-empty piercing sequences; the corresponding function $J_0(\tau)$ is shown in Figure \ref{357nonspin} on the left.

The only difference with the self-conjugate case is that in this situation the plot is not symmetric under $\tau\mapsto 1-\tau$, so that the two intervals, which are
\begin{equation*}
[0.0710,0.0728]\text{ and }[0.3557,0.3575]
\end{equation*}
have to be studied independently. On the other hand, in all non self-conjugate examples we have worked out, the situation turns out to be symmetric with respect to another point in $\mathbb{T}$; in the example under consideration the symmetry point is
\begin{equation*}
0.5-2/7\approx 0.2143.    
\end{equation*}
This is to be expected because all the manifolds in the Hodgson-Weeks census arise as branched double covers of links in $S^3$, and given that in our examples $b_1(Y)=1$, the covering involution acts on $H^1(Y;\mathbb{Z})=\mathbb{Z}$ as multiplication by $-1$.
\par
\begin{figure}
\centering
  \includegraphics[width=.6\linewidth]{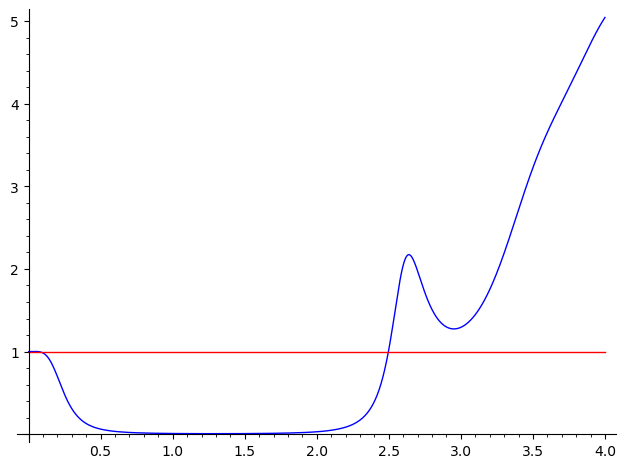}
  \caption{The plot of $J_s(\tau_1)$ for the pair of self conjugate spin$^c$ structures under consideration on \#357, computed at cutoff $R=7$.}
  \label{357nonspinbooker}
\end{figure}

The plot of $J_s(\tau_1)$, where $\tau_1$ is the midpoint of the first interval, is shown in Figure \ref{357nonspinbooker}, while the plot of $\Gamma_K(\tau)$ is shown in Figure \ref{357nonspin} on the right. The same exact approach of the previous sections allows us to conclude that there is exactly a single transverse crossing in each interval.

\vspace{0.3cm}
\subsection{Small spectral gap in \#3250. }
\begin{figure}
\centering
\begin{subfigure}{.45\linewidth}
  \includegraphics[width=\linewidth]{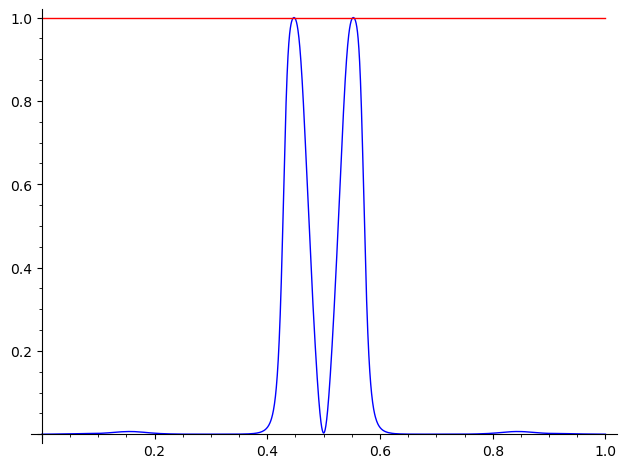}

\end{subfigure}
\begin{subfigure}{.45\linewidth}
  \includegraphics[width=\linewidth]{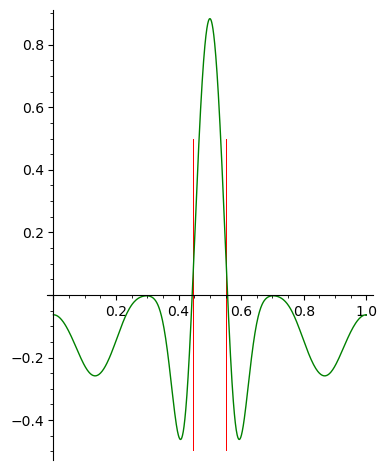}
\end{subfigure}
\caption{For one of the (self-conjugate) spin$^c$ structures on \#3250, the plot of $J_0(\tau)$ computed at cutoff $R=7.5$, and the plot of $\Gamma_{K_7}(\tau)$.}
 \label{3250first}
\end{figure}

\begin{figure}
\centering
\begin{subfigure}{.45\linewidth}
  \includegraphics[width=\linewidth]{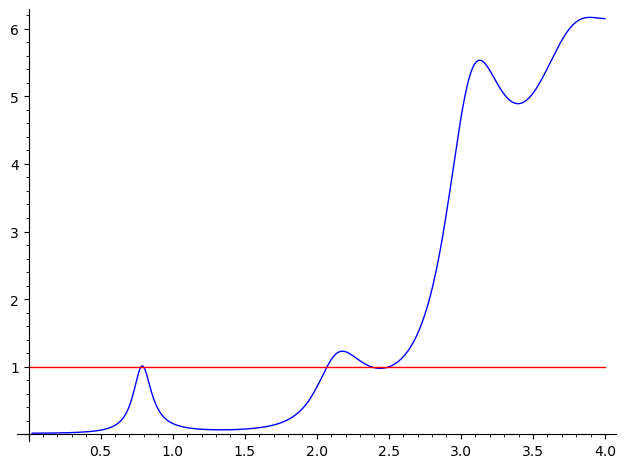}
\end{subfigure}
\begin{subfigure}{.45\linewidth}
  \includegraphics[width=\linewidth]{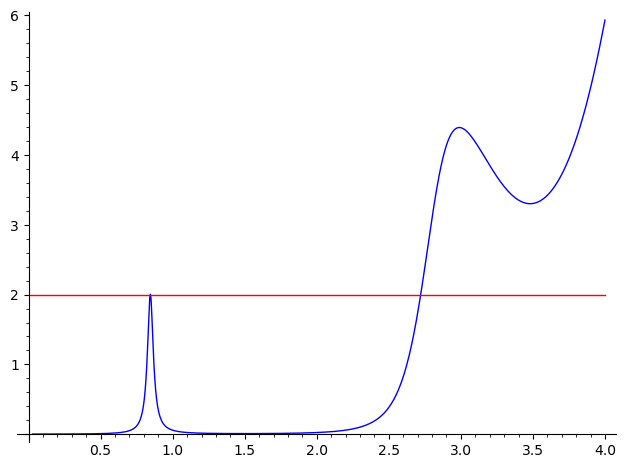}
\end{subfigure}
  \caption{The plots of $J_s(0.4)$ and $J_s(0.5)$ for the spin$^c$ structure on \#3250 under consideration, computed at cutoff $R=7.5$.}
  \label{3250crossing}
\end{figure}

We next consider one of the two spin$^c$ structures on the manifold \#3250. The plots of $J_0(\tau)$, computed for $R=7.5$ is shown in Figure \ref{3250first} in the left, and one shows that in the intervals
\begin{equation*}
[0.4467,0.4480]\cup[0.5520,0.5533]
\end{equation*}
there is exactly one small eigenvalue as in Section \ref{smalleigint}; the plot of the family of odd trace formula (again computed with $K_7$) is given in Figure \ref{3250first} on the right. The higher peak of $\Gamma_{K_7}$ is due to the fact that at $\tau=0.5$ we are considering a genuine spin (rather than spin$^c$) connection, so that the corresponding Dirac operator is quaternionic and eigenvalue multiplicities are always even (see also Figure \ref{3250crossing} on the right).

\vspace{0.3cm}

Again, to check that the crossing actually occurs in the interval, we look at nearby values; see in Figure \ref{3250crossing} the plots of $J_s(0.4)$ and $J_s(0.5)$ respectively. The main new feature that may cause worries for our estimates is that for $J_s(0.4)$, the lower bound on the next eigenvalues is not as good as in our previous examples; because of this, we check that
\begin{equation*}
\frac{1}{2}\sum_{i\in\mathbb{Z}}(\sinc^7)'(s_i(0.4))=\Gamma_{K_7}(0.4)\approx -0.4511,
\end{equation*}
while, using that $|s_i(0.4)|\geq 2$ for $i\neq0$, we bound

\begin{equation*}
\left|\sum_{i\neq0}(\sinc^7)'(s_i(0.4))\right|\leq 0.1911,
\end{equation*}
and conclude again that $s_0(0.4)>0$.

\vspace{0.3cm}

\begin{figure}
\centering
  \includegraphics[width=0.6\linewidth]{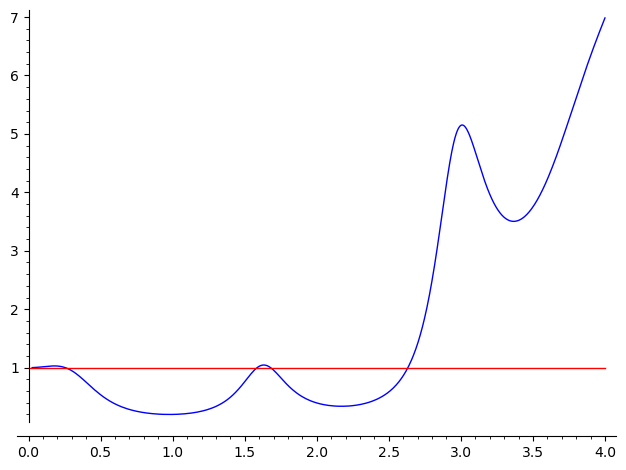}
  \caption{The plot of $J_s(0.44735)$ for the spin$^c$ structure on \#3250 under consideration, computed at cutoff $R=7.5$.}
  \label{3250smalleig}
\end{figure}

The main complication in the proof of transversality in this example is that near the crossing the next eigenvalue is not as large as in the case of Census 357 we discussed in detail in Section \ref{singlecrossing}. In fact, the plot of $J_s(\tau_0)$ at the midpoint $\tau_0=0.44735$ suggests that the manifold admits a eigenvalue with absolute value in the interval $[1.5856,1.6799]$, see Figure \ref{3250smalleig}. Indeed, the methods of Subsection \ref{upperbounddensity} (using the fact that $\sinc^6(x)\geq 0.9977$ for $|x|\leq 0.04715$, the latter being half the length of the interval) prove that there is at most one eigenvalue in this interval.

\vspace{0.3cm}
We can then employ the same strategy of proof as before to prove that 
\begin{equation*}
s_0'(\tau)\neq 0\text{ for }\tau\in [0.4467,0.4480].
\end{equation*}
In particular, focusing for simplicity on the midpoint $\tau_0$, we compute
\begin{equation*}
\widetilde{\Gamma}_{K_7}(\tau_0)\approx 23.1598.
\end{equation*}
We then bound, using the fact that the next eigenvalue is $\geq 2.5$ in absolute value,
\begin{equation*}
\left\lvert \sum_{i\neq 0}\widehat{K_7}'(s(\tau_0))\right\rvert \leq 0.6424+0.0675=0.7099,
\end{equation*}
where the first summand, the maximum of $(\sinc^7)''$ in $[1.5856,1.6799]$, comes from the second smallest eigenvalue, while the second one comes from the tail of the sum. For \#3250 we have $C_Y\leq 8.8899$ (cf. Appendix \ref{C0bounds}), so that $$|s_i'(\tau_0)|\leq 2\pi C_Y\leq 55.8569,$$ and therefore
\begin{equation*}
\left\lvert\sum_{j\neq 0} \widehat{K_7}'(s_i(t_0))s_i'(t_0)\right\rvert\leq 0.7099\times 55.8569\leq 39.6529,
\end{equation*}
which is strictly less that $2|\widetilde{\Gamma}_{K_7}(0.44735)|$, and we conclude.

\vspace{0.3cm}
\subsection{Spin$^c$ structure with four crossings on \#10867. }The example in Theorem \ref{thm2} has significantly larger volume, and to determine the piercing sequences we computed the length spectrum up to $R=8.5$, which took roughly a week.
\par
\begin{figure}
\centering
\begin{subfigure}{.45\linewidth}
  \includegraphics[width=\linewidth]{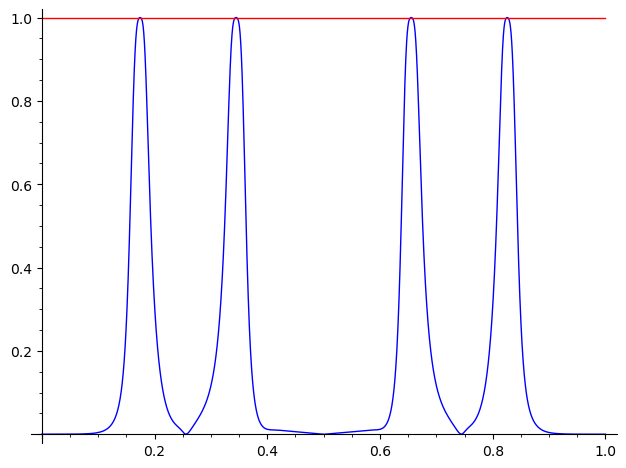}
\end{subfigure}
  \begin{subfigure}{.45\linewidth}
  \includegraphics[width=\linewidth]{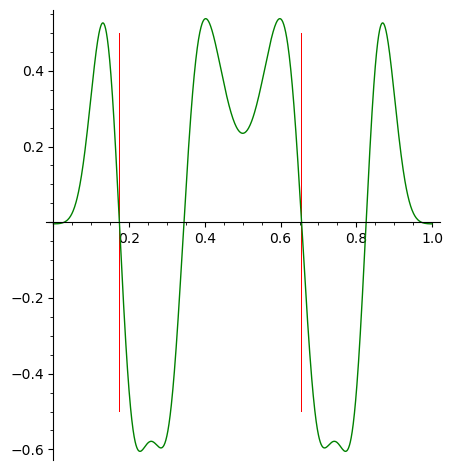}
  \end{subfigure}
   \caption{The plots of $J_0(\tau)$ and $\Gamma_K(\tau)$ for one of the self-conjugate spin$^c$ structures on \#10867, computed using $R=8.5$.}
   \label{10867book}
\end{figure}

The corresponding function $J_0(\tau)$ (evaluated using the whole range) and the plot of $\Gamma_{K}(t)$ where
\begin{equation}\label{stretchedtest}
K=K_8\left(\frac{x}{1.0625}\right)\text{ with }
K_8(x)=x\cdot \left(\frac{1}{2}\mathbf{1}_{[-1,1]}(x)\right)^{*8},
\end{equation}
i.e. $K$ is $K_8$ stretched to have support in $[-8.5,8.5]$, can be found in Figure \ref{10867book}. From this, the intervals with small eigenvalues are determined to be
\begin{equation*}
I_1=[0.1735,0.1751],\quad I_2=[0.3438,0.3455]    
\end{equation*}
and the corresponding ones $I_3, I_4$ under the symmetry $\tau\mapsto 1-\tau$.

\vspace{0.3cm}

\begin{figure}
\centering
\begin{subfigure}{.45\linewidth}
  \includegraphics[width=\linewidth]{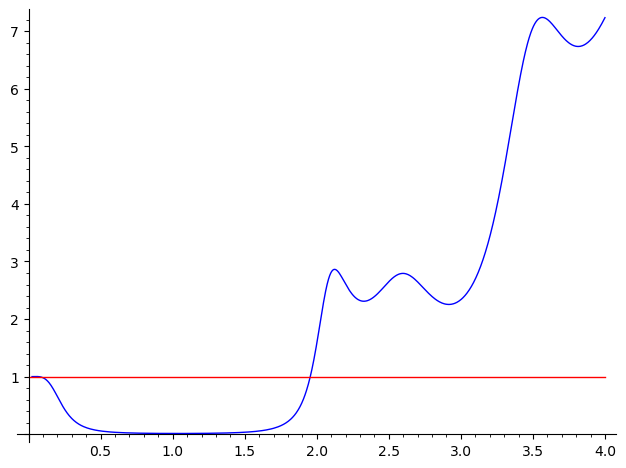}
\end{subfigure}
  \begin{subfigure}{.45\linewidth}
  \includegraphics[width=\linewidth]{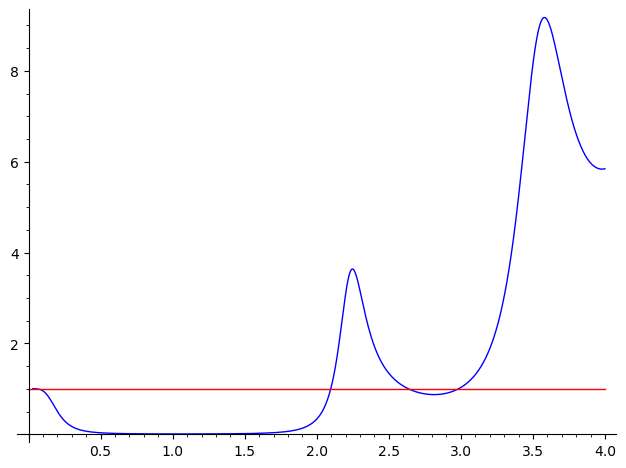}
  \end{subfigure}
   \caption{The plots of $J_s(\tau)$ for the spin$^c$ structure under consideration on \# 10867 at the midpoints of the intervals $I_1$ and $I_2$, computed at cutoff $R=8.5$.}
   \label{10867bookermidpoints}
\end{figure}

The plot of the function $J_s(\tau)$ at the midpoints of the intervals are shown in Figure \ref{10867bookermidpoints}. In the case of the intervals $I_2$ and $I_3$, there is a good spectral gap and few small eigenvalues, and one proves easily that there is a single transverse crossing. The case of $I_1$ and $I_4$ is significantly more challenging because of the smaller spectral gap and the higher eigenvalue density, and we will work it out in detail (focusing on $I_1$). First of all, one has
\begin{equation*}
J_{1.9541}(\tau)\leq 0.9983 \text{ for }\tau\in I_1,
\end{equation*}
so that the second smallest eigenvalue always has absolute value $\geq1.9541$. Applying the method of subsection \ref{upperbounddensity} to bound spectral density in intervals (using the eight, rather than sixth, convolution power) one furthermore shows that for all $\tau\in I_1$, the number of spectral parameters in the short intervals
\begin{equation}\label{intervalbound}
[1.9541,2.2941],\quad [2.2941,2.4741],\quad [2.4741,2.7041]
\end{equation}
is at most $4$, $5$, and $6$ respectively.
\par
In our convention for the Fourier transform we have
\begin{align*}
\widehat{K}(t)&=1.0625\cdot \widehat{K_8}(1.0625\cdot t)\\
\widehat{K}'(t)&=1.0625^2\cdot \widehat{K_8}'(1.0625\cdot t)
\end{align*}
so that
\begin{equation*}
-i\widehat{K}'(t)=1.0625^2\cdot (\sinc^8)''(1.0625\cdot t).
\end{equation*}
We have that for the geometric side of the derivative of the odd trace formula
\begin{equation*}
\widetilde{\Gamma}_K(\tau)\in[-19.9527,-19.9260]\text{ for all }\tau\in I_1.
\end{equation*}
On the other hand we can bound
\begin{equation*}
\left|\sum_{i\neq 0}\widehat{K}'(s_i(t))\right|\leq 4\times 0.06735+5\times 0.00402+6\times 0.00045+0.0023\leq 0.2945,
\end{equation*}
where the first three terms are bounds for the contribution of the eigenvalues in the short intervals (\ref{intervalbound}), while the last term bounds the contribution of the eigenvalues $\geq 2.7041$ computed as in the previous examples (using the convolution eight power to bound spectral densities). We conclude the proof of Theorem \ref{thm2} because for \#10867 we have $C_Y\leq  5.7836$ (cf. Appendix \ref{C0bounds}).

\vspace{0.5cm}

\appendix
\section{Explicit upper bounds on the $C^0$ norm of closed $1$-forms}\label{C0bounds}
In Section \ref{traceformulas} we exploited the universal bound of Lemma \ref{eigenbound} to justify taking the derivative of the family of the trace formulas. On the other hand, in order to perform the tail estimates required in the proofs of our main results, we need \textit{explicit} upper bounds on the value of $C_Y$; the goal of this section is to discuss an algorithm and its implementation that constructs closed $1$-forms representing a generator $H^1(Y;\mathbb{Z})=\mathbb{Z}$ with explicitly controlled norm; the output can be found in Table 2.

\begin{table}[ht]
\centering
\begin{tabular}{|c| c |}
\hline
Label & upper bound on $C_Y$ \\ [0.5ex] 
\hline
356&3.9921\\
357&3.5151\\
381&9.5088\\
735&3.2401\\
882&4.8550\\
1155& 5.6521\\
1280&9.0190\\
1284&2.6791\\
3250&8.8899\\ 
3673& 4.4709\\
10867& 5.7836\\
\hline
\end{tabular}
\caption{Upper bounds on $C_Y$; we are omitting \#734 and \#790 because there are no crossings to be certified.}
\end{table}

\begin{remark}
In the opposite direction, there are two natural ways to provide \textit{lower} bounds on $C_Y$ in terms of the geometry and topology of $Y$. First, if $\gamma$ is a closed geodesic, denoting $[\gamma]\in H_1(Y;\mathbb{Z})/\mathrm{tors}=\mathbb{Z})$ we have
\begin{equation*}
|[\gamma]|=\left|\int_{\gamma}\alpha\right|\leq C_Y\cdot \ell(\gamma)
\end{equation*}
so that
\begin{equation*}
C_Y\geq\sup\frac{[\gamma]}{{\ell(\gamma)}},
\end{equation*}
the supremum being taken over all closed geodesics. Similarly, taking a circle-valued primitive $f:Y\rightarrow \mathbb{R}/\mathbb{Z}$ of $\alpha$, the coarea formula implies that
\begin{equation*}
 \int_{\mathbb{R}/\mathbb{Z}}\mathrm{Area}(f^{-1}(t))dt\leq \int_Y|\alpha|\leq C_Y\cdot\vol(Y).
\end{equation*}
With some care about the singularities of $f$, the fact that in for a class in $H_2(Y;\mathbb{Z})$ the least area norm is at least $\pi$ times the Thurston norm \cite{BD} implies then that
\begin{equation*}
C_Y\geq \frac{\pi\cdot\mathrm{Th}(Y)}{\vol(Y)},
\end{equation*}
where we denote by $\mathrm{Th}(Y)$ the Thurston norm of the generator of $H_2(Y;\mathbb{Z})$.
\end{remark}
In fact, our method will work more in general for $b_1>0$ and any class $a\in H^1(Y;\mathbb{Z})$. We take as input a Dirichlet domain $\mathsf{D}\subset\mathbb{H}^3$ for $Y$ together with the face-pairing maps $\{g_i\}_{i\in I}$, which are elements in $\mathrm{Isom}^+(\mathbb{H}^3)$ (as computed for example by SnapPy). First of all, we interpret the class $a$ as a homomorphism
\begin{equation*}
\varphi_a: \pi_1(Y)\rightarrow \mathbb{Z}.
\end{equation*}
Our concrete goal is then to describe an explicit function
\begin{equation}\label{aclass}
F:\mathsf{D}\rightarrow \mathbb{R}
\end{equation}
with the property that
\begin{equation}\label{equivariant}
F(g_i\cdot x)=F(x)+\varphi_a(g_i)\text{ for all }x\in\partial\mathsf{D}\text{ and } i\in I,
\end{equation}
so that the closed $1$-form $dF$ represents $a$.
\par
In order to do this, we first we will fix a triangulation $\mathcal{T}=\{T_i\}$ of $\mathsf{D}$ in simplices with geodesic boundary, and consider functions obtained by assigning values at the vertices of the triangulation, and extended in a canonical `linear' fashion.
\\
\par
In what follows, we first discuss the geometry of such linear extensions, and then discuss the specific details of our implementation. An important simplifying assumption we will make is that all the tetrahedra in our triangulation admit a \textit{circumcenter}, i.e. the four vertices are equidistant from some point in $\mathbb{H}^3$; unlike the euclidean case, this is a non-trivial condition.

\vspace{0.5cm}
\textbf{Geometry of linear extensions in hyperbolic space.} A geodesic tetrahedron $T$ in hyperbolic space comes with a natural parametrization
\begin{equation*}
\eta:\Delta\rightarrow T
\end{equation*}
from the standard $3$-simplex
\begin{equation*}
\Delta=\left\{(t_0,t_1,t_2,t_3)\lvert t_i\geq 0 \text{ and }\sum t_i\leq 1\right\}\subset\mathbb{R}^4
\end{equation*}
given by convex combinations, the definition of which we now recall (see for example \cite{Mar}). We consider the hyperboloid model $\mathsf{I}$ for hyperbolic space. To set notation, we consider the standard Minkowski space given by $\mathbb{R}^4$ with coordinates $x=(x_0,x_1,x_2,x_3)$ and Lorentzian inner product
\begin{equation*}
\langle x,y\rangle=-x_0y_0+\sum_{i=1}^3x_iy_i.
\end{equation*}
Then $\mathsf{I}$ is the upper half of the two-sheeted hyperboloid
\begin{equation*}
\mathsf{I}=\left\{x\lvert \langle x,x\rangle=-1, x_0>0\right\}
\end{equation*}
equipped with restriction of the Lorentzian inner product, which is a Riemannian metric of constant curvature $-1$. In this model, the geodesics are given exactly by the intersection of $\mathsf{I}$ with planes $\mathbb{R}^4$, and the distance between points $x,y\in \mathsf{I}$ is given by
\begin{equation*}
d(x,y)=\cosh^{-1}(-\langle x, y\rangle).
\end{equation*}
Denoting
\begin{equation*}
\|x\|=\sqrt{-\langle x,x\rangle},
\end{equation*}
we have the hyperbolic analogue of the radial projection map
\begin{align*}
P: \{x\lvert  \langle x,x\rangle<0,x_0>0\}&\rightarrow \mathsf{I}\\
x&\mapsto \frac{x}{\|x\|}.
\end{align*}
Given two distinct points $x,y\in \mathsf{I}$, the geodesic segment in $\mathsf{I}$ between them is parametrized (not at constant speed) by
\begin{equation}\label{convex}
P((1-t)x+ty) \text{ for }t\in[0,1].
\end{equation}
In particular, given four points $v_i$ in $\mathsf{I}$ (which are linearly independent in $\mathbb{R}^4$), the geodesic simplex $T$ in $\mathsf{I}$ they determine is parametrized via the convex combinations
\begin{align*}
\eta: \Delta&\rightarrow T\\
(t_0,t_1,t_2,t_3)&\mapsto P(\sum t_iv_i).
\end{align*}
Of course, this is the composition of the projection $P$ with the linear parametrization of the affine simplex spanned by the $v_i$ in $\mathbb{R}^4$. Notice that this parametrization is equivariant with respect to the action of $\mathrm{Isom}^+(\mathbb{H}^3)=\mathrm{SO}_0(3,1)$. This is because a point $p\in \Delta$ is the projection of a unique point in the convex hull in $\mathbb{R}^4$ of the vertices of $\Delta$, and isometries act via linear maps on such convex hull.
\\
\par
Suppose now we are given a function $f$ on the vertices of the geodesic tetrahedron $T$. This naturally extends to a function
\begin{equation}\label{linext}
F: T\rightarrow \mathbb{R}
\end{equation}
given as the composition
\begin{equation}\label{compositionLip}
T\stackrel{\eta^{-1}}{\longrightarrow}\Delta{\longrightarrow}\mathbb{R}
\end{equation}
where the latter is the unique affine linear function on $\Delta$ agreeing with $f\circ\eta$ on the vertices of $\Delta$; we will refer to $F$ as the linear extension of $f$. The goal of this section is then to estimate from above the norm of $dF$, or equivalently the Lipschitz constant of $F$, in terms of $f$ and $T$, under the simplifying assumption that $T$ admits a circumcenter; this can be achieved with the following steps.
\\
\par
\textit{Step 0: find the circumcenter of $T$. }This is done via linear algebra by noticing that the circumcenter $c$ of $T$ lies at the intersection of the bisectors of pairs of vertices $v_i,v_j$, which are described by the equations
\begin{equation*}
\langle z, v_i\rangle=\langle z,v_j\rangle.
\end{equation*}
Notice that the existence of the circumcenter is equivalent to the fact that the one-dimensional space of solutions to these equations in $\mathbb{R}^4$ intersects the hyperboloid $\mathsf{I}$.
\\
\par
\textit{Step 1: move the circumcenter to the origin $o=(1,0,0,0)$.}
We can apply isometry that maps $c$ to $o$ (e.g. one can just use the reflection at the midpoint of the segment). Under this isometry, the the tetrahedron $T$ is mapped to $\tilde{T}$, a geodesic tetrahedron with circumcenter at the origin $o$. In particular, its vertices $\tilde{v_i}$ lie on the boundary of the ball $B_R\subset\mathsf{I}$ of radius $R$ around $o$; by equivariance, the natural parametrization
\begin{equation*}
\tilde{\eta}: \Delta\rightarrow \tilde{T}
\end{equation*}
is obtained by composing the isometry with $\eta$. Because of this, we can assume from now on that our vertices lie on the boundary of $B_R$.
\\
\par
\textit{Step 2: fixing a Euclidean structure on $\Delta$.} In order to evaluate the Lipschitz constant of (\ref{linext}), we will factor it as in (\ref{compositionLip}) and use the chain rule. To do so, we \textit{identify} $\Delta$ with the affine span of the $\{v_i\}$, equipped with the Riemannian metric induced by the standard Euclidean norm $x_0^2+x_1^2+x_2^2+x_3^2$ on $\mathbb{R}^4$, which we denote as $\Delta_{\mathrm{hor}}$, where the subscript reminds us that $\Delta_{\mathrm{hor}}$ is horizontal (i.e. $x_0$ is constant). We then have to study the Lipschitz constants of the maps 
\begin{equation}\label{compositionLip}
T\stackrel{P^{-1}_{\mathrm{hor}}}{\longrightarrow}\Delta_{\mathrm{hor}}{\longrightarrow}\mathbb{R}
\end{equation}
where to simplify the notation we denote the inverse of $P:\Delta_{\mathrm{hor}}\rightarrow T$ by $P^{-1}_{\mathrm{hor}}$. Of course, the Lipschitz constant of the second map is directly computed in terms of linear algebra, so the remaining non-trivial computations consists of providing an upper bound on the Lipschitz constant of $P^{-1}_{\mathrm{hor}}$.
\\
\par
\begin{figure}
  \centering
\def\svgwidth{0.8\textwidth}
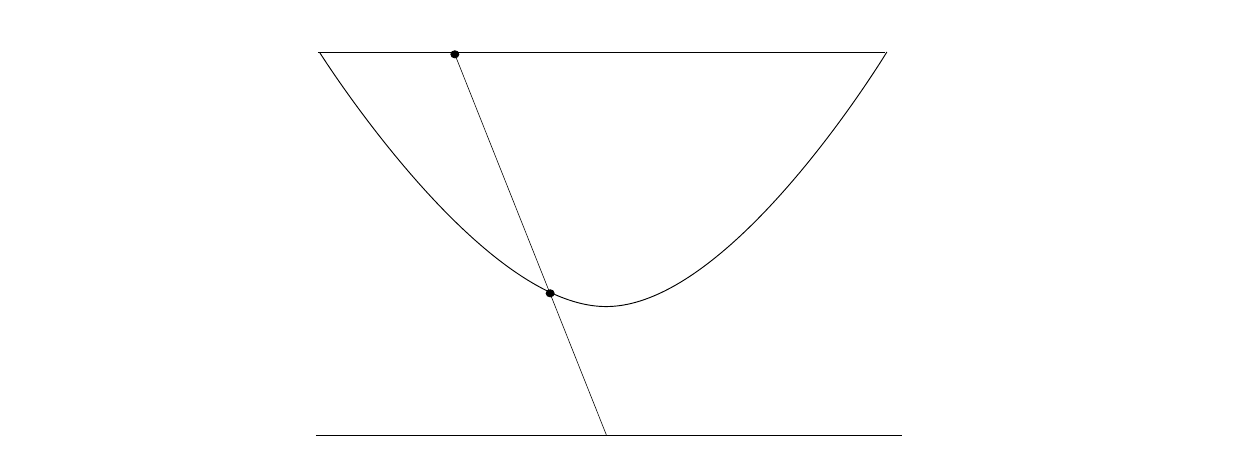
\caption{A schematic picture of the map $P_R^{-1}$.}
\label{PR}
\end{figure}

\textit{Step 3: the Lipschitz constant of $P_{\mathrm{hor}}^{-1}$.}
We show that the Lipschitz constant of $P_{\mathrm{hor}}^{-1}$ is $\cosh(R)$ by considering more in general the map $P_R^{-1}$ from the ball of radius $R$ around $o$ in $\mathsf{I}$ to the horizontal disk with the same boundary (see Figure \ref{PR}) as follows.
\par
Let us introduce coordinates $(z,t,\vartheta)$ on $\mathbb{R}^{4}$ where $z=x_0$ and $(t,\vartheta)$ are polar coordinates in the $\mathbb{R}^{3}$ corresponding to the last $3$ coordinates. Consider the points
\begin{equation*}
(\sqrt{1+a^2}, a,\vartheta)\text{ for }\vartheta\in S^2
\end{equation*}
on $\mathsf{I}$, which form the boundary of the ball $B_R$ of radius
\begin{equation*}
R=\cosh^{-1}(\sqrt{1+a^2})
\end{equation*}
around the origin $o=(1,0,0,0)$. The horizontal ball in this situation is just the disk
\begin{equation*}
(\sqrt{1+a^2},t,\vartheta)\text{ for }t\in[0,a]\text{ and }\vartheta\in S^2,
\end{equation*}
equipped with the standard metric Euclidean
\begin{equation*}
dt^2+t^2 g_{S^{2}}.
\end{equation*}
In these coordinates projection map $P_R$ from the horizontal ball to $B_R$ is given by
\begin{equation*}
(t,\vartheta)\mapsto \left(\frac{\sqrt{1+a^2}}{\sqrt{1+a^2-t^2}}, \frac{t}{\sqrt{1+a^2-t^2}},\vartheta \right).  
\end{equation*}
The exponential map at the origin in hyperbolic space is given by
\begin{equation*}
(s,\vartheta)\mapsto (\sinh(s),\vartheta,\cosh(s));
\end{equation*}
recall that in these coordinates the hyperbolic metric is given by
\begin{equation*}
g_{\mathrm{hyp}}=ds^2+\sinh^2(s)g_{S^{n-1}}.
\end{equation*}
In particular, denoting by $s(t)$ the $s$-coordinate of $P(t,\vartheta)$, we have
\begin{equation*}
s(t)=\mathrm{arcsinh}\left(\frac{t}{\sqrt{1+a^2-t^2}}\right),
\end{equation*}
so that the pullback of the hyperbolic metric to the horizontal disk is given by
\begin{equation*}
P_R^*(g_{\mathrm{hyp}})=\frac{a^2+1}{(1+a^2-t^2)^2}dt^2+\frac{t^2}{1+a^2-t^2}g_{S^{n-1}}.
\end{equation*}
From this we see that the differential of the inverse $P_R^{-1}$ multiplies lengths by
\begin{equation*}
\frac{(1+a^2-t^2)}{\sqrt{a^2+1}}\text{ and }\sqrt{1+a^2-t^2}
\end{equation*}
in the radial and tangential directions respectively, and both of these quantities have maximum $\sqrt{1+a^2}=\cosh(R)$ at the origin.

\vspace{0.5cm}

\textbf{Implementation.} The Dirichlet domains constructed by SnapPy for our manifolds of interest are all relatively small, and we can obtain a geodesic triangulation by simply considering the one given by barycentric subdivision of faces and edges and coning on the basepoint of the domain. This does not always lead to a triangulation with all tetrahedra having circumcenter, but for all our examples we find such triangulations by changing the basepoint of the Dirichlet domain. This leads to triangulations for our manifolds with around 150-200 geodesic tetrahedra, all of which have circumcenter.
\par
We then assign to the vertices of the triangulation random values satisfying all the conditions (\ref{equivariant}), and extend linearly to the simplices as in the previous section to obtain a continuous function on $Y$. The $C^0$-norm of $dF$ on each simplex can be bounded above explicitly (say by $C>0$) by the discussion of the previous section (the circumradii of the tetrahedra are quite small for our examples, so the distortion between the euclidean and hyperbolic norms is in fact not too large). Notice though that in general $dF$ is only piecewise continuous; on the other hand, for any $\varepsilon>0$ one readily constructs by interpolation along vertices, edges and faces a function $\tilde{F}$ which is $C^1$ and for which $$\|d\tilde{F}\|_{C^0}\leq C+\varepsilon,$$so that we can simply consider the norm of $dF$ for our purposes as we round up the bound at the very end of the process.
\par
The constant $C$ already provides one explicit upper bound on the quantity $C_Y$ we are interested in. To obtain a (much) better upper bound, we optimize within the space of choices, i.e. assignment of values to the vertices. Unfortunately, this problem cannot be solved via linear programming techniques because the computation of Lipschitz constants of the latter functions on tetrahedra involve significant non-linearities. The approach we take to optimizing: select a random line in the space of choices passing through $F$, minimize the target quantity on this line to obtain an improved function $F^\ast$, then replace $F$ by $F^\ast$ and iterate.  The computations in Table 2 were obtained after roughly 5000-15000 iterations of the procedure.

\vspace{0.5cm}

\bibliography{biblio.bib}
\bibliographystyle{alpha}
\end{document}

%% file: tetrahedra.pdf_tex
\begingroup%
  \makeatletter%
  \providecommand\color[2][]{%
    \errmessage{(Inkscape) Color is used for the text in Inkscape, but the package 'color.sty' is not loaded}%
    \renewcommand\color[2][]{}%
  }%
  \providecommand\transparent[1]{%
    \errmessage{(Inkscape) Transparency is used (non-zero) for the text in Inkscape, but the package 'transparent.sty' is not loaded}%
    \renewcommand\transparent[1]{}%
  }%
  \providecommand\rotatebox[2]{#2}%
  \newcommand*\fsize{\dimexpr\f@size pt\relax}%
  \newcommand*\lineheight[1]{\fontsize{\fsize}{#1\fsize}\selectfont}%
  \ifx\svgwidth\undefined%
    \setlength{\unitlength}{595.27559055bp}%
    \ifx\svgscale\undefined%
      \relax%
    \else%
      \setlength{\unitlength}{\unitlength * \real{\svgscale}}%
    \fi%
  \else%
    \setlength{\unitlength}{\svgwidth}%
  \fi%
  \global\let\svgwidth\undefined%
  \global\let\svgscale\undefined%
  \makeatother%
  \begin{picture}(1,0.38095238)%
    \lineheight{1}%
    \setlength\tabcolsep{0pt}%
    \put(0,0){\includegraphics[width=\unitlength,page=1]{tetrahedra.pdf}}%
    \put(0.45095338,0.14798477){\color[rgb]{0,0,0}\makebox(0,0)[lt]{\lineheight{1.25}\smash{\begin{tabular}[t]{l}$x$\end{tabular}}}}%
    \put(0.35132414,0.35357715){\color[rgb]{0,0,0}\makebox(0,0)[lt]{\lineheight{1.25}\smash{\begin{tabular}[t]{l}$P_{R}^{-1}(x)$\end{tabular}}}}%
  \end{picture}%
\endgroup%

%% file: piercing.bbl
\begin{thebibliography}{CDGW}

\bibitem[APS76]{APS3}
M.~F. Atiyah, V.~K. Patodi, and I.~M. Singer.
\newblock Spectral asymmetry and {R}iemannian geometry. {III}.
\newblock {\em Math. Proc. Cambridge Philos. Soc.}, 79(1):71--99, 1976.

\bibitem[BD17]{BD}
Jeffrey~F. Brock and Nathan~M. Dunfield.
\newblock Norms on the cohomology of hyperbolic 3-manifolds.
\newblock {\em Invent. Math.}, 210(2):531--558, 2017.

\bibitem[BS07]{BS}
Andrew~R. Booker and Andreas Str\"ombergsson.
\newblock Numerical computations with the trace formula and the {S}elberg eigenvalue conjecture.
\newblock {\em J. Reine Angew. Math.}, 607:113--161, 2007.

\bibitem[Bus10]{Bus}
Peter Buser.
\newblock {\em Geometry and spectra of compact {R}iemann surfaces}.
\newblock Modern Birkh\"auser Classics. Birkh\"auser Boston, Ltd., Boston, MA, 2010.
\newblock Reprint of the 1992 edition.

\bibitem[CDGW]{SnapPy}
Marc Culler, Nathan Dunfield, Matthias Goerner, and Jeff Weeks.
\newblock Snappy, a computer program for studying the geometry and topology of 3-manifolds.

\bibitem[Cha84]{Cha}
Isaac Chavel.
\newblock {\em Eigenvalues in {R}iemannian geometry}, volume 115 of {\em Pure and Applied Mathematics}.
\newblock Academic Press, Inc., Orlando, FL, 1984.
\newblock Including a chapter by Burton Randol, With an appendix by Jozef Dodziuk.

\bibitem[Che23]{Chen}
Jacopo Chen.
\newblock Computing the twisted $l^2$-euler characteristic.
\newblock {\em to appear in Groups Geom. Dyn.}, 2023.

\bibitem[HW]{Census}
Craig Hodgson and Jeff Weeks.
\newblock A census of closed hyperbolic 3-manifolds.

\bibitem[Kat95]{Kato}
Tosio Kato.
\newblock {\em Perturbation theory for linear operators}.
\newblock Classics in Mathematics. Springer-Verlag, Berlin, 1995.
\newblock Reprint of the 1980 edition.

\bibitem[KM07]{KM}
Peter Kronheimer and Tomasz Mrowka.
\newblock {\em Monopoles and three-manifolds}, volume~10 of {\em New Mathematical Monographs}.
\newblock Cambridge University Press, Cambridge, 2007.

\bibitem[Lin16]{LinLec}
Francesco Lin.
\newblock Lectures on monopole {F}loer homology.
\newblock In {\em Proceedings of the {G}\"okova {G}eometry-{T}opology {C}onference 2015}, pages 39--80. G\"okova Geometry/Topology Conference (GGT), G\"okova, 2016.

\bibitem[Lin24a]{LinTheta}
Francesco Lin.
\newblock Monopole {F}loer homology and invariant theta characteristics.
\newblock {\em J. Lond. Math. Soc. (2)}, 109(5):Paper No. e12895, 27, 2024.

\bibitem[Lin24b]{LinDir}
Francesco Lin.
\newblock Topology of the {D}irac equation on spectrally large three-manifolds.
\newblock {\em preprint}, 2024.

\bibitem[LL21]{LL3}
Francesco Lin and Michael Lipnowski.
\newblock Closed geodesics and fr\o yshov invariants of hyperbolic three-manifolds.
\newblock {\em to appear in JEMS}, 2021.

\bibitem[LL22a]{LL2}
Francesco Lin and Michael Lipnowski.
\newblock Monopole {F}loer homology, eigenform multiplicities, and the {S}eifert-{W}eber dodecahedral space.
\newblock {\em Int. Math. Res. Not. IMRN}, (9):6540--6560, 2022.

\bibitem[LL22b]{LL1}
Francesco Lin and Michael Lipnowski.
\newblock The {S}eiberg-{W}itten equations and the length spectrum of hyperbolic three-manifolds.
\newblock {\em J. Amer. Math. Soc.}, 35(1):233--293, 2022.

\bibitem[Mar]{Mar}
Bruno Martelli.
\newblock An introduction to {G}eometric {T}opology.

\bibitem[Roe98]{Roe}
John Roe.
\newblock {\em Elliptic operators, topology and asymptotic methods}, volume 395 of {\em Pitman Research Notes in Mathematics Series}.
\newblock Longman, Harlow, second edition, 1998.

\bibitem[SN]{SN}
J.~J. Sakurai and Jim Napolitano.
\newblock Modern quantum mechanics.

\bibitem[Tur98]{Tur}
Vladimir Turaev.
\newblock A combinatorial formulation for the {S}eiberg-{W}itten invariants of {$3$}-manifolds.
\newblock {\em Math. Res. Lett.}, 5(5):583--598, 1998.

\end{thebibliography}
